\documentclass[12pt]{article}
\usepackage{amssymb}
\newtheorem{theorem}{Theorem}[section]
\newtheorem{lemma}[theorem]{Lemma}
\newtheorem{proposition}[theorem]{Proposition}

\newcommand{\qed}{\hfill $\square$\vskip .2cm}

\newcommand{\sect}[1]{\section{#1}\setcounter{equation}{0}}
\newcommand{\comment}[1]{}
\newcommand{\Real}{\mbox{\rm I\kern-.23em\hbox{R}}}

\newcommand{\be}{\begin{equation}}
\newcommand{\ee}{\end{equation}}

\def\v{\vskip .5cm}
\def\Z{\mathbb Z}
\def\F{\mathfrak{F}}
\def\R{{I\!\!R}}

\def\B{{\mathcal B}}

\def\A{{\mathcal{A}}}
\def\S{{\mathcal{S}}}

\def\L{{\mathcal{L}}}
\def\N{{\mathcal{N}}}
\def\l{\mathfrak{L}}
\def\s{\mathfrak{S}}
\def\a{\mathfrak{A}}
\def\<{\langle}
\def\>{\rangle}
\def\Epsilon{{\mathcal E}}
\def\f{{\mathfrak f}}
\def\g{{\mathfrak g}}

\def\G{{\mathcal G}}
\def\C{{\mathcal C}}
\def\T{{\mathcal T}}
\def\n{{\mathfrak N}}

\begin{document}

\author{Sunder Sethuraman}

\title{Diffusive variance for a tagged
particle in $d\leq 2$ asymmetric simple exclusion} 

\thispagestyle{empty} \maketitle \abstract{The study of equilibrium
fluctuations of a tagged particle in
    finite-range simple exclusion processes has a long history. The
    belief is that the scaled centered tagged particle motion behaves
    as some sort of homogenized random walk. In fact, invariance
    principles have been proved in all dimensions $d\geq 1$ when the
    single particle jump rate is unbiased, in $d\geq 3$ when the jump
    rate is biased, and in $d=1$ when the jump rate is in addition
    nearest-neighbor.
    
    The purpose of this article is to give some partial results in the
    open cases in $d\leq 2$.  Namely, we show the tagged particle
    motion is ``diffusive'' in the sense that upper and lower bounds
    are given for the tagged particle variance at time $t$ on order
    $O(t)$ in $d=2$ when the jump rate is biased, and also in $d=1$
    when in addition the jump rate is not nearest-neighbor.
    Also, a characterization of the
 tagged particle variance is given. 
 The main methods are in analyzing $H_{-1}$ norm variational
 inequalities.}
\v
  
\thanks{$*$ Research partially supported by NSF/DMS-0504193 and NSA-H982300510041} \\

{\sl AMS (2000) subject classifications}: Primary 60K35; secondary
 60F05.\\
 {\sl Address and Email}: 396 Carver Hall, Department of Mathematics, Iowa State
 University, Ames, IA \ 50011; sethuram@iastate.edu
\eject

\sect{Introduction and Results} One of the interesting questions in
Spitzer's seminal paper on particle systems \cite{Spitzer} asks for
the asymptotics of a distinguished or ``tagged,'' particle as it
interacts with others. Although the tagged particle is not in
general Markovian, due to the particle interactions, the
understanding is that it behaves in some sense as a ``homogenized''
random walk. In the context of finite-range translation-invariant
simple exclusion processes, this belief has been substantiated in
large part through a quilt of results sometimes depending on the
specific form of the single particle jump rate $p$, and the
dimension $d$ of the underlying lattice $\Z^d$.

For instance, laws of large numbers, both in equilibrium
\cite{Saada-se} and non-equilibrium \cite{Reza-lln} have been shown.
Also, equilibrium central limit theorems and invariance principles
when $p$ is mean-zero \cite{Arratia},\cite{RV}, \cite{KV},
\cite{Varadhan}, and when $p$ has a drift in $d\geq 3$ \cite{SVY}
and in $d=1$ when $p$ is in addition nearest-neighbor \cite{Kipnis}
have been proved. See also \cite{LOV},\cite{LV} for fluctuations in
$d=1$ with respect to a non-translation invariant $p$.  Non-trivial
non-equilibrium fluctuation results have even been derived in $d\geq
1$ when $p$ is symmetric (excluding the $d=1$ nearest-neighbor case)
\cite{Reza-pr}, and recently in the exceptional case in $d=1$ when
$p$ is symmetric and nearest-neighbor \cite{JL}.  In addition, large
deviations results have been proved in some cases \cite{QRV},
\cite{Sepa}. Some of these results and others are reviewed in
\cite{F}, section 4.VIII \cite{Liggett}, chapter 4.III
\cite{Liggett2}, chapter 6 \cite{Spohn}, and sections 4.3, 8.4 and
11.5 \cite{KL}.

In terms of equilibrium fluctuations, however, open are the
behaviors in $d= 2$ when $p$ has a drift, and also in $d=1$ when in
addition $p$ is not nearest-neighbor.  The difficulty in their
solution is roughly that in low dimensions with asymmetry one
has to deal with more involved particle interactions than in high
dimensions, where transience estimates can be used, and under symmetry,
when reversibility helps.  The main goal of this article is to shed
light on the open low dimensional cases by giving some upper and lower bounds on the
variance of the tagged particle at time $t$ which are ``diffusive,''
that is on order $O(t)$ (Theorems \ref{maintheorem1} and
\ref{maintheorem2}).
In addition, a characterization of the variance, which recasts an
expression in the literature (cf. equation (1.18) \cite{DeMF}) in
terms of certain ``dynamical'' and ``static'' contributions, is
given (Theorem \ref{variance_characterization}).

The method of the upper bounds is to bound above the variance of a
``drift'' additive functional as $O(t)$ by estimating certain
$H_{-1}$ variational formulas with the help of integral estimates in
the spirit of Bernardin's work for occupation times \cite{B}.  In
particular, one of the main contributions of this article is to give
a framework for tagged particle $H_{-1}$ norms in which
``environment'' and ``tagged-shift'' dynamics are understood. The
variance characterization, and lower bounds follow from explicit
computations, and comparisons with ``symmetrized'' variances as in
Loulakis \cite{Loulakis}.

Loosely speaking, the simple exclusion process follows the motion of
a collection of random walks on the lattice $\Z^d$ in which jumps to
already occupied vertices are suppressed.  More precisely, let
$\Sigma = \{0,1\}^{\Z^d}$ and let $\eta(t) \in \Sigma$ represent the
state of the process at time $t$. That is, the configuration at time
$t$ is given in terms of occupation variables $\eta(t) =
\{\eta_i(t): i\in \Z^d\}$ where $\eta_i(t)=0 {\rm \ or \ }1$
according to whether the vertex $i\in \Z^d$ is empty or full at time
$t$.  Let $p=\{p(i,j): i,j\in \Z^d\}$ be the single particle
transition rates. Throughout this article we concentrate on the
translation-invariant finite-range case: $p(i,j)=p(0,j-i)=p(j-i)$
and $p(x)=0$ for $|x|>R$ and an integer $R<\infty$.  In addition, to
avoid technicalities, we concentrate on the situation when
$(p(i)+p(-i))/2$ is irreducible, and $p(0)=0$.  We will say $p$ is
nearest-neighbor when the range $R=1$.

The system $\eta(t)$ is a Markov process on $D(\R_+,\Sigma)$ with
semi-group $T_t$ and generator, well defined on functions $\phi$
supported on a finite number of vertices, namely ``local''
functions, \be \label{generator} (L\phi)(\eta) \ =\  \sum_{i,j\in \Z^d}
\eta_i(1-\eta_j)p(j-i)(\phi(\eta^{i,j})-\phi(\eta)) \ee where
$\eta^{i,j}$ is the ``exchanged'' configuration, $(\eta^{i,j})_i
=\eta_j$, $(\eta^{i,j})_j = \eta_i$ and $(\eta^{i,j})_k = \eta_k$
for $k\neq i,j$. We note the transition rate
$\eta_i(1-\eta_j)p(j-i)$ for $\eta \rightarrow \eta^{i,j}$
represents the exclusion property.

With respect to a configuration $\eta$, distinguish now one of the
particles and call it the tagged particle.  Let $x(t)\in \Z^d$ be
its position at time $t$. To compensate for the non-Markovian
character of the tagged motion, we form the larger process
$(x(t),\eta(t))$ which is Markovian. In fact, as is standard
practice, we will consider the system in the reference frame of the
tagged particle, $(x(t),\zeta(t))$ where
$\zeta(t)=\pi_{x(t)}\eta(t)$.  Here, for a configuration $\eta\in
\Sigma$, the $k$-shifted state is $\pi_k\eta$ where $(\pi_k\eta)_l =
\eta_{k+l}$ for $l\in \Z^d$. The ``reference frame'' process
$\zeta(t)$ is also Markovian with semi-group $\T_t$, and generator
$\L$ well defined on local functions,
\begin{eqnarray*}
(\L\phi)(\zeta) &=& \sum_{i,j\in \Z^d\setminus\{0\}}
\zeta_i(1-\zeta_j)p(j-i)
(\phi(\zeta^{i,j})-\phi(\zeta))\\
&&\ \ \ \ \ \ \ + \sum_{j\in \Z^d\setminus\{0\}}
(1-\zeta_j)p(j)(\phi(\tau_{j}\zeta))-\phi(\zeta))
\end{eqnarray*}
where $\tau_j \zeta = \pi_j (\zeta^{0,j})$ accounts for the
reference frame shift when the tagged particle displaces by $j$.

Naturally, $\L$ splits as $\L = \L^e +\L^t$ where $(\L^e\phi)(\zeta)
= \sum_{i,j\in \Z^d\setminus\{0\}} \zeta_i(1-\zeta_j)p(j-i)
(\phi(\zeta^{i,j})-\phi(\zeta))$ and $(\L^t\phi)(\zeta) = \sum_{j\in
\Z^d\setminus\{0\}}
(1-\zeta_j)p(j)(\phi(\tau_{j}\zeta))-\phi(\zeta))$ correspond to
movement around, and by the tagged particle, e.g. ``environment''
and ``tagged-shift'' motions, respectively. The main idea of the
reference process is that, although the tagged particle is always at
the origin ($\zeta_0(t) \equiv 1$), one can keep track of the
position of the tagged particle by counting the various reference
``$j$-shifts'' (cf. (\ref{martdecomp})). We refer to \cite{Liggett}
for details of the construction of these processes.

We now discuss the equilibria for these systems. Let $P_\rho$, for
$\rho \in [0,1]$, be the infinite Bernoulli product measure over
$\Z^d$ with coin-tossing marginal $P_\rho\{\eta_i
=1\}=1-P_\rho\{\eta_i=0\}=\rho$. It is known that $P_\rho$ and
$Q_\rho =P_\rho(\cdot|\zeta_0=1)$ are invariant extremal measures
for $L$ and $\L$ respectively \cite{Saada-se}. We remark with
respect to $P_\rho$, the semi-group $T_t$ and generator $L$ can be
extended to $L^2(P_\rho)$ (cf. section IV.4 \cite{Liggett});
similarly, with respect to $Q_\rho$, $\T_t$ and $\L$ can be extended
to $L^2(Q_\rho)$. We note the adjoints $L^*$ and $\L^*$ with respect
to $P_\rho$ and $Q_\rho$, corresponding to time-reversal, are
straightforwardly computed and identified as generators
corresponding to reversed jump rates $p(-\cdot)$. It will sometimes
be convenient to write $\L$ into symmetric and anti-symmetric parts,
$\L = \S + \A$ where $\S = (\L + \L^*)/2$ and $\A = (\L -\L^*)/2$.
We note the operator $\S$ is the generator of a reference frame
process with symmetric jump rates $(p(\cdot)+p(-\cdot))/2$.  Also,
as before, $\S$ and $\A$ can be split into ``environment'' and
``tagged-shift'' parts, $\S = \S^e +\S^t$ and $\A=\A^e +\A^t$.

We denote $E_\rho$ for expectation with respect to the reference
process measure starting from $Q_\rho$.  Denote also, for
vector-valued functions $f,g:\Sigma \rightarrow \R^m$ and $m\geq 1$,
the innerproduct $\<f,g\>_\rho = E_\rho[f\cdot g]$, and $L^2$ norm
$\|f\|_0= \sqrt{\<f,f\>_\rho}$ with respect to $Q_\rho$.

We now specify a family of martingales associated with the exclusion
process.  For $j\in \Z^d$, let $N_{j}(t)$ denote the counting
processes which count the number of $j$-shifts made by the reference
process, e.g. $j$-displacements of the tagged particle, up to time
$t\geq 0$. By subtracting appropriate compensators, we can then form
the martingale $M_j(t) = N_j(t) - A_j(t)$ where $A_j(t) = \int_0^t
p(j)(1-\zeta_j(s))ds$.  These martingales, as jumps are not
simulateneous, are orthogonal for $j\in \Z^d$. Then, the tagged
particle position $x(t)$ may be written into the sum of a martingale
and an additive functional term, $$x(t) \ =\ \sum_j jN_j(t) \ = \
\sum_j jM_j(t) + \sum_j jA_j(t).$$ These relations, by stationarity
of the process measure, give the quadratic variation $E_{\rho}
[M^2_j(t)]=(1-\rho)p(j)t$ and mean position, $E_{\rho}[x(t)]=
(\sum_j jp(j))(1-\rho)t$. Then, after centering,
\begin{equation}
\label{martdecomp}x(t) - E_\rho[x(t)] \ = \  M(t) +
A(t)\end{equation} with martingale $M(t) = \sum jM_j(t)$ and
``drift'' $A(t)=\int_0^t \F(\zeta(s))ds$ with $\F(\zeta) = \sum
jp(j)(\rho-\zeta_j)$.



Let now
$$V(t) \ = \ E_{\rho}\bigg[|x(t)-E_{\rho}[x(t)]|^2\bigg].$$
Define also the
measure
 $d\mu_{k,\rho} =
 (\zeta_{k}/\rho)dQ_\rho$ and its expectation $E_{k,\rho}$ for $k\in \Z^d\setminus\{0\}$.
 The
first result is a characterization of the variance.  In a different
form, it was first derived by De Masi and Ferrari (cf. equation
(1.18) \cite{DeMF}), however, the interpretation below seems new.
See also \cite{Szrtg} for analogous expressions in zero-range
processes.

 \begin{theorem}
 \label{variance_characterization}
 In $d\geq 1$,
  $$
 V(t) \ = \
 (1-\rho)\sum_j |j|^2 p(j) t  +2\rho\sum_j  jp(j) \cdot \int_0^t
 \bigg\{
 E_\rho[x(s)] - E_{-j,\rho}[x(s)]\bigg\}ds .$$
 \end{theorem}

 The first term above, $(1-\rho)\sum |j|^2 p(j)
 t$, is the mean quadratic variation of the martingale $M(t)$ and
can be thought of as a ``dynamical'' part of the variation. The
second term, however, as a difference in expected tagged
 particle positions from different initial measures, is in a
 sense variation due to initial conditions.

 We note in $d=1$ when $p$
 is totally asymmetric and nearest-neighbor, say $p(1)>0$ and $p(i)=0$
 for $i\neq 1$, the second term in the decomposition vanishes as
 the ``extra'' particle at $-1$, being behind, cannot interfere with the tagged
 particle position; in this case, $V(t) = p(1)(1-\rho)t$ and moreover it is known the tagged motion
 is actually a Poisson process with rate $p(1)(1-\rho)$ (cf. Corollary VIII.4.9 \cite{Liggett}).
   Also, in $d=1$ when $p$ is nearest-neighbor, the formula can be
 evaluated to some extent, and the limit $\lim_{t\rightarrow \infty} V(t)/t = (1-\rho)|p(1)-p(-1)|$ has been proved \cite{DeMF}.

However, for the next upper
 bounds, other methods are used.
\begin{theorem}
\label{maintheorem1} When $p$ has a drift, $\sum jp(j)\neq 0$, in
$d=2$, and in $d=1$ when additionally $p$ is not nearest-neighbor,
we have a constant $C=C(d,p,\rho)$ such that
 $$ V(t)\ \leq \ C t.$$
\end{theorem}

For a general lower bound, we only give an estimate on a
``Tauberian'' quantity which resembles $V(t)$.
\begin{theorem}\label{maintheorem2}
In $d\geq 1$ and for $0\leq \rho <1$, excluding the nearest-neighbor
symmetric case in $d=1$ when $p(1)=p(-1)$, we have a constant
$C=C(d,p,\rho)>0$ such that
$$\liminf_{\lambda \downarrow 0} \ \lambda^2 \int_0^\infty e^{-\lambda t}V(t)dt \ \geq \ C.$$
\end{theorem}

The lower bound, by formal (non-rigorous) analogies, suggests
$$\frac{1}{T} \int_0^\infty e^{-t/T}V(t)dt \ \sim \ \frac{1}{T}
\int_0^{T}V(t)dt \
 \sim \ V(T) \ \geq \ C T.$$
We note also our proofs of Theorems \ref{maintheorem1} and
\ref{maintheorem2} only give gross estimates on the constants
$C(d,p,\rho)$.

 However, well-known when $p$
 is mean-zero and not nearest-neighbor in $d=1$,
 biased in $d\geq 3$, or biased and nearest-neighbor in $d=1$ , the variance is on
 order $V(t)=O(t)$
 \cite{KV},\cite{Varadhan}, \cite{SVY}, \cite{Kipnis}; in
the excluded $d=1$ nearest-neighbor symmetric case, due to
``trapping'' phenomena, $V(t) =O(\sqrt{t})$ \cite{Arratia}. Also,
when $\rho = 1$, there is no motion and $V(t)\equiv 0$.

%


We remark now, in terms of remaining open questions,
the limit
\begin{equation} \label{V(t)limit} \lim_{t\rightarrow \infty}
\frac{1}{t}V(t) \ = \ \sigma^2(d,p,\rho), \end{equation} and full
invariance principles should hold more generally in $d\leq 2$ when
$\sum jp(j)\neq 0$.

We suspect more detailed $H_{-1}$ norm estimation might allow
martingale approximation of the tagged position $x(t)$ leading to
limits (\ref{V(t)limit}) and invariance principles in this
situation.  Namely, one wants to show the ``drift'' $\F$ (cf.
(\ref{martdecomp})) can be approximated in terms of $\L u_\epsilon$
where $u_\epsilon$ is a local function satisfying $\|\F - \L
u_\epsilon\|_{H_{-1}}<\epsilon$. This type of program was done in
\cite{SVY} in $d\geq 3$ using ``transience estimates'' which
unfortunately are not available in $d\leq 2$. We hope however the
basic $H_{-1}$ estimates given in this article will serve as
building blocks for subsequent work.

The structure of the article is to prove first the variance
characterization and lower bound in section 2.  The upper bound is
proved in section 4 with the aid of some preliminaries in section 3
and technical computations in section 5.

\sect{Proofs of Theorems \ref{variance_characterization} and \ref{maintheorem2}}
Let $s$ and $a$ be the symmetric and anti-symmetric parts of $p$,
$s(i)=(p(i)+p(-i))/2$ and $a(i)=(p(i)-p(-i))/2$ for $i\in \Z^d$.
Recall the ``drift'' function $\F= \sum jp(j)(\rho - \zeta_j)$ in
the introduction, and define analogous ``drifts'' $\F_s(\zeta) =
\sum js(j)(\rho-\zeta_j)$ and $\F_{\leftarrow}(\zeta)= \sum
jp(-j)(\rho-\zeta_j)$ corresponding to rates $s(\cdot)$ and
$p(-\cdot) = s(\cdot)-a(\cdot)$ respectively.\vskip .2cm

{\it Proof of Theorem \ref{variance_characterization}.}
Following
decomposition (\ref{martdecomp}), write
\begin{eqnarray}
\label{prop2.2line} V(t) &=&(1-\rho)\sum_j |j|^2p(j)t +
2E_\rho[M(t)\cdot A(t)] +
E_\rho[|A(t)|^2]\nonumber\\
&=&(1-\rho)\sum_j |j|^2p(j)t + 2\int_0^t E_\rho[M(s)\cdot
\F(\zeta(s))]ds +
E_\rho[|A(t)|^2]\nonumber\\
&=&(1-\rho)\sum_j |j|^2p(j)t + 2\int_0^t E_\rho[x(s)\cdot
\F(\zeta(s))]ds
\end{eqnarray}
where we note $E_\rho[|A(t)|^2] = 2\int_0^tE_\rho[A(s)\cdot
\F(\zeta(s))]ds$. We now reverse time at $s$, and note the
time-reversed process $\zeta(s-\cdot)$ with respect to process
measure started from $Q_\rho$ has the same distribution as the
process with reversed jump rates. In particular, $N_j(s)$ with
respect to the process begun from $Q_\rho$ has the same distribution
as $N_{-j}(t)$ with respect to the reversed process.  Hence, as
$x(t) =\sum jN_j(t)$, we have $E_\rho[x(s)\cdot \F(\zeta(s))] =
E^*_\rho[-x(s)\cdot \F(\zeta(0))]$ where $E^*_\rho$ is expectation
with respect to the reversed process begun with $Q_\rho$. Then, by
spatial reflection, simple manipulations, and recalling the measure
$d\mu_{k,\rho}= (\zeta_k/\rho)dQ_\rho$ with expectation
$E_{k,\rho}$, we have
\begin{eqnarray} \label{expression} -E^*_\rho[x(s)\cdot
\F(\zeta(0))] & = &
-E_\rho\bigg[\sum_k kN_{-k}(s) \cdot \sum_j jp(j)(\rho-\zeta_{-j}(0))\bigg]\\
 & = & \rho\sum_j
jp(j)\cdot \bigg\{E_\rho[x(s)] -
E_{-j,\rho}[x(s)]\bigg\}.\nonumber\end{eqnarray}\qed

\vskip .1cm

{\it Proof of Theorem \ref{maintheorem2}.}  The proof follows
straightforwardly from Propositions \ref{symlowerbound} and
\ref{lowerboundprop} below which allow comparisons with the tagged
particle variance for the symmetrized process. \qed

  Let $E^\S_\rho$
be expectation with respect to the symmetric reference process
generated by $\S$ with initial distribution $Q_\rho$.  Let also
$V_s(t) = E^\S_\rho[|x(t) - E_\rho[x(t)]|^2]$ be the corresponding
variance of the tagged particle at time $t$. Then, the following
estimate is proved in \cite{KV}.

\begin{proposition}
\label{symlowerbound} In $d\geq 1$ and for $0\leq \rho <1$, except
for the nearest-neighbor symmetric case in $d=1$ when $p(1)=p(-1)$,
we have a constant $C=C(d,p,\rho)>0$ such that $V_s(t) \geq Ct$ for
all $t\geq 0$.
\end{proposition}


Form now, for $\lambda >0$, two resolvent equations,
$$\lambda u_\lambda -\L u_\lambda = \F \ \ \ \ {\rm and \ \ \ \ } \lambda
v_\lambda - \S v_\lambda = \F_s$$ with respect to $u_\lambda =
(\lambda -\L)^{-1}\F = \int_0^\infty e^{-\lambda t}(T_t\F) dt$ and
$v_\lambda = (\lambda -\S)^{-1}\F_s$. We now state a comparison, in
whose proof, the last part is Corollary 1 \cite{Loulakis}.
\begin{proposition}
\label{lowerboundprop} We have
\begin{eqnarray*}
\int_0^\infty e^{-\lambda t}[V(t) - V_s(t)]dt & = &
\frac{2}{\lambda^2}\bigg[ \<\F_s,(\lambda - \S)^{-1}\F_s\>_\rho
-\<\F_{\leftarrow},(\lambda -\L)^{-1}\F\>_\rho\bigg]\\
&=&\frac{2}{\lambda^2}\bigg[\lambda \|u_\lambda -v_\lambda \|_0^2 +
\<u_\lambda-v_\lambda, (-\S)(u_\lambda -v_\lambda)\>_\rho\bigg
].\end{eqnarray*}\end{proposition}
We note, as $-\S$ is a
non-negative operator, the Dirichlet form $\<u_\lambda-v_\lambda,
(-\S)(u_\lambda -v_\lambda)\>_\rho \geq 0$, and so as a consequence,
 $\int_0^\infty e^{-\lambda t}V(t)dt \geq
\int_0^\infty e^{-\lambda t}V_s(t)dt$. \vskip .2cm

 {\it Proof.} We first evaluate further (\ref{expression}) as
 $$-E_\rho[x(s)\cdot \F_\leftarrow(0)] \ = \ -E_\rho[A(s)\cdot
\F_\leftarrow(0)]$$ after the martingale part in $x(s)=M(s)+A(s)$
vanishes. Then, the last term of (\ref{prop2.2line}) equals
$$
 - 2\int_0^t E_\rho[A(s)\cdot
\F_{\leftarrow}(\zeta(0))]ds\\
\ =\  - 2\int_0^t \int_0^s \<\F_\leftarrow, T_s\F\>_\rho drds.$$
Hence, by two integration by parts,
\begin{eqnarray*}
\int_0^\infty e^{-\lambda t}V(t)dt &=& \lambda^{-2}(1-\rho)\sum_j
|j|^2p(j)- 2\lambda^{-2}\int_0^\infty e^{-\lambda t}
\<\F_{\leftarrow},T_t \F\>_\rho dt\\
&=& \lambda^{-2}(1-\rho)\sum_j |j|^2p(j)
 -2\lambda^{-2}\<\F_{\leftarrow},(\lambda
-\L)^{-1}\F\>_\rho. \end{eqnarray*}
 Since, $\sum |j|^2 p(j) = \sum
|j|^2s(j)$ and $\F = \F_{\leftarrow}=\F_s$ when $p(\cdot)=s(\cdot)$, we obtain the first equality in the proposition
directly.

For the second equality, we compute, using $\F+\F_{\leftarrow} =
2\F_s$, the two resolvent equations and $\<u_\lambda, \A
u_\lambda\>_\rho = 0$, that
\begin{eqnarray*}
\<\F_s, v_\lambda\>_\rho - \<\F_{\leftarrow},u_\lambda\>_\rho
&=&\<\F_s, v_\lambda\>_\rho + \<\F,u_\lambda\>_\rho -
2\<\F_s,u_\lambda\>_\rho\\
&=&\<v_\lambda,(-\S)v_\lambda\>_\rho +\<u_\lambda,(-\S)u_\lambda\>_\rho\\
&&\ \ \ \ \ \ \ \ \ \ \ \ \ \ \ \ \   +\lambda\|u_\lambda\|_0^2
+\lambda \|v_\lambda\|^2_0 - 2\<\
\F_s,u_\lambda\>_\rho\\
&=&\<v_\lambda,(-\S)v_\lambda\>_\rho +\<u_\lambda,(-\S)u_\lambda\>_\rho\\
&&\ \ \ \ \ \ \ \ \ \ \ \ \ \ \ \ \  +\lambda \|u_\lambda
-v_\lambda\|_0^2 +2\lambda\<v_\lambda ,u_\lambda\>_\rho -
2\<\F_s,u_\lambda\>_\rho.
\end{eqnarray*}
Since $2\<\lambda v_\lambda,u_\lambda\>_\rho -
2\<\F_s,u_\lambda\>_\rho = -2\<(-S)v_\lambda,u_\lambda \>_\rho$,
 we have the
right-side equals $\lambda \|u_\lambda -v_\lambda \|_0^2 +
\<u_\lambda-v_\lambda, (-\S)(u_\lambda -v_\lambda)\>_\rho$ as
desired. \qed

\sect{Preliminaries for Upper Bound}
We discuss here some definitions and results useful for the
upperbound.

\subsection{Duality} As the tagged particle is
always at the origin with respect to the reference process, consider
the underlying lattice $\Z^d\setminus\{0\}$ .  Let $\Epsilon_d$
denote the collection of finite subsets of $\Z^d\setminus\{0\}$, and
let $\Epsilon_{d,n}$ be those subsets of cardinality $n\geq 0$. Let
$\beta_\rho = \sqrt{\rho(1-\rho)}$ and, for non-empty $B\in
\Epsilon_d$, let $\Psi_B$ be the function
$$\Psi_B(\zeta) \ =\ \prod_{x\in B} \frac{\zeta_x -
  \rho}{\beta_\rho}$$
when $0<\rho<1$, and $\Psi_B \equiv 0$ when $\rho =0$ or $1$.
 By convention, we set $\Psi_\emptyset \equiv 1$. One can check that
  $\{\Psi_B: B\in \Epsilon_d\}$ is a Hilbert basis of $L^2(Q_\rho)$.  In
  particular, any function $f\in L^2(Q_\rho)$ has decomposition
$$f \ = \ \sum_{n\geq 0} \sum_{B\in \Epsilon_{d,n}} \f(B)\Psi_B$$
with coefficient $\f:\Epsilon_d \rightarrow \R$ which in general
depends on $\rho$. Then, for $f,g\in L^2(Q_\rho)$, we define
innerproduct
$$\<\f,\g\>\ = \ \<f,g\>_\rho \ = \ \sum_{B\in \Epsilon_d}\f(B)\g(B),$$
and $L^2$ norm by $\|\f\|_0^2 = \|f\|_0^2 = \<f,f\>_\rho$.

Let also $\C_{d,n}$ be the subspace of coefficient functions on
$\Epsilon_{d,n}$. When $f$ is in the span of $\{\Psi_B: B\in
\Epsilon_{d,n}\}$, we have $\f\in \C_{d,n}$, and we say both $f$ and
its coefficient $\f$ are of degree $n$. Note also, when $f$ is
local, then $\f$ is also local on $\Epsilon_d$, that is with support
on a finite number of subsets of $\Z^d\setminus\{0\}$.

The operators $\L$, $\S$ and $\A$ have counterparts $\l=\l^e +\l^t$,
$\s=\s^e +\s^t$ and $\a=\a^e +\a^t$ which act on ``coefficient''
functions $\f$:
$$\L^e f \ = \ \sum_{B\in \Epsilon} (\l^e\f)(B)\Psi_B, \ \ \S^e f \ = \
\sum_{B\in \Epsilon} (\s^e\f)(B)\Psi_B, \ \ {\rm and \ \ } \A^e f \
= \ \sum_{B\in \Epsilon} (\a^e\f)(B)\Psi_B$$ with analogous
expressions for $\l^t$, $\s^t$ and $\a^t$.

Recall the symmetric and anti-symmetric parts of $p$, $s(i)= (p(i)
+p(-i))/2$ and $a(i) = (p(i)-p(-i))/2$ for $i\in \Z^d$; by
assumption $s(0)=a(0)=0$.  For $B\subset \Z^d\setminus\{0\}$, denote
$$B_{x,y} = \left\{\begin{array}{rl}
B\setminus \{x\} \cup \{y\} & \ {\rm when \ } x\in B, y\not\in B\\
B\setminus \{y\} \cup \{x\} & \ {\rm when \ } x\not\in B, y\in B\\
B&\ {\rm otherwise}\end{array}\right.
$$
and $$\tau_x B \ =\ \left\{\begin{array}{rl} B +x & \ {\rm when \ }
-x\not\in B\\
(B+x)\setminus \{0\}\cup \{x\}&\ {\rm when \ } -x\in
B\end{array}\right. $$ where as usual $B+x = \{i+x: i\in B\}$ for
$B$ nonempty, and $\emptyset +x =\emptyset$. As in \cite{SVY}, the
symmetric parts $\s^e$ and $\s^t$ can be computed as
\begin{eqnarray*}
(\s^e\f)(B) & = & \frac{1}{2}\sum_{x,y\in \Z^d\setminus\{0\}}
s(y-x)[\f(B_{x,y}) -
\f(B)]\\
(\s^t\f)(B) &=& (1-\rho)\sum_{z\not\in
B\atop z\in \Z^d\setminus\{0\}}s(z)[\f(\tau_{-z}B)-f(B)] + \rho\sum_{z\in B}s(z)[f(\tau_{-z}B)-f(B)]\\
&& \ + \beta_\rho\sum_{z\not\in B\atop z\in \Z^d\setminus\{0\}}
s(z)[\f(B\cup \{z\}) -
\f(\tau_{-z}(B\cup \{z\}))]\\
&&\ + \beta_\rho\sum_{z\in B} s(z)[\f(B\setminus \{z\}) -
\f(\tau_{-z}(B\setminus\{z\}))].
\end{eqnarray*}
 Note that $\s^e\f\in \C_{d,n}$ for $\f\in \C_{d,n}$, and so $\s^e$
``preserves'' degrees. However, $\s^t$ does not ``preserve'' degrees
but, as will be seen, we will not need to deal directly with $\s^t$
in our calculations.

Also, the anti-symmetric parts $\a^e$ and $\a^t$ are decomposed into
sums of three operators which preserve, increase, and decrease the
degree of the function acted upon: $\a^e = \a^e_0 + \a^e_+ +\a^e_-$
and $\a^t = \a^t_0 +\a^t_+ +\a^t_-$ where
\begin{eqnarray*}
(\a^e_0\f)(B) & =& (1-2\rho)\sum_{x\in B\atop y\not\in B, y\in
\Z^d\setminus\{0\}}
a(y-x)[\f(B_{x,y})-\f(B)]\\
(\a^e_+\f)(B) & = & 2\beta_\rho\sum_{x,y\in B} a(y-x) \f(B\setminus\{y\})\\
(\a_-^e\f)(B) &= & -2\beta_\rho\sum_{x,y\not\in B\atop x,y\in
\Z^d\setminus\{0\}}
a(y-x) \f(B\cup \{x\})\\
(\a^t_0\f)(B) & = & (1-\rho)\sum_{z\not\in B\atop z\in
\Z^d\setminus\{0\}}a(z)[\f(\tau_{-z}B)-
\f(B)] +\rho \sum_{z\in B}a(z)[\f(\tau_{-z}B) - \f(B)]\\
(\a^t_+\f)(B) &= & \beta_\rho\sum_{z\in
B}a(z)[\f(B\setminus\{z\}) - \f(\tau_{-z}(B\setminus\{z\}))]\\
(\a^t_-\f)(B) &=& \beta_\rho\sum_{z\not\in B\atop z\in
\Z^d\setminus\{0\}}a(z)[\f(B\cup\{z\}) - \f(\tau_{-z}(B\cup\{z\}))].
\end{eqnarray*}
It will also be helpful to write $\a$ in terms of its explicit
``degree'' actions,
$$\a \ = \ \sum_{n\geq 0} \bigg ( \a_{n,n-1} + \a_{n,n} +
\a_{n,n+1}\bigg )$$ where $\a_{m,n}$ is the part which takes a
degree $m$ function to a degree $n$ function. Here, by convention
$\a_{0,-1} \equiv 0$ is the zero operator; one also sees $\a_{0,0} =
\a_{1,0} = \a_{0,1} \equiv 0$.
Similarly, $\a^e$ and $\a^t$ can be decomposed in terms of degree
actions $\a^e_{m,n}$ and $\a^t_{m,n}$ so that $\a_{m,n} = \a^e_{m,n}
+\a^t_{m,n}$ for $m,n\geq 0$.  We later evaluate in Proposition
\ref{FT}, and its proof in section 5, some of the relevant actions.

\subsection{Variational Formulas} 
Define, for $\lambda >0$ and local $\phi$, the $H_{1,\lambda,\L}$
norm $\|\cdot\|_{1,\lambda,\L}$ by
$$\|\phi\|_{1,\lambda,\L}^2 \ = \ \<\phi, (\lambda -\S)\phi\>_\rho + \<\A\phi, (\lambda -
\S)^{-1}\A\phi\>_\rho$$ where we note $\<\phi,(-\S)\phi\>_\rho,
\<\A\phi,(\lambda -\S)^{-1}\A\phi\>_\rho\geq 0$ as $-\S$ is a
non-negative operator. The $H_{1,\lambda,\L}$ Hilbert space is then
the completion over local functions with respect to this norm.

To define a dual norm, consider for $f\in L^2(Q_\rho)$ and local
$\phi$ that $$\<f,\phi\>_\rho \ \leq \ \|f\|_0\|\phi\|_0 \ \leq\
\lambda^{-1/2}\|f\|_0\|\phi\|_{1,\lambda, \L} \ .$$  Then, the dual
norm of $\|\cdot\|_{1,\lambda, \L}$, given by
$$
\|f\|_{-1,\lambda, \L} \ = \ \sup_{\phi {\rm \ local}\atop
\|\phi\|_{1,\lambda, \L} = 1} \<f,\phi\>_\rho, $$ is always finite
with bound $\|f\|^2_{-1,\lambda, \L} \leq \lambda^{-1} \|f\|_0^2$.
Let $H_{-1,\lambda, \L}$ be the corresponding Hilbert space with
respect to $\|\cdot\|_{-1,\lambda, \L}$. An equivalent expression
for $\|f\|_{-1,\lambda, \L}$, given in the next result, is proved in
p. 46-47 \cite{Olla}.

\begin{proposition}
\label{infvarexpression}
For $f\in L^2(Q_\rho)$ and $\lambda >0$, we have
\begin{eqnarray*}
\|f\|^2_{-1,\lambda,L}
&=& \<f,(\lambda
-\L)^{-1}f\>_\rho \\
&=&\sup_{g\ {\rm local}}\bigg\{ 2\< f,g\>_\rho - \<g,(\lambda
-\S )g\>_\rho -\<\A g ,(\lambda -\S )^{-1}\A g \>_\rho\bigg\}\\
 &=&\inf_{g\ {\rm local}}\bigg\{\< f - \A g, (\lambda -\S)^{-1}(f-\A g)\>_\rho +
\< g,(\lambda -\S)g\>_\rho\bigg\}.
\end{eqnarray*}
\end{proposition}

Hence, when $\L =\S$ is symmetric, we have for $f$ local,
$\|f\|^2_{1,\lambda,\S} = \<f,(\lambda -\S)f\>_\rho$ and
$\|f\|^2_{-1,\lambda,\S} = \<f,(\lambda-\S)^{-1}f\>_\rho$.  In this
context, it will be useful to define corresponding $H_{1}$ and
$H_{-1}$ ``coefficient'' norms, that is, $\|\f\|^2_{1,\lambda,\s} =
\<\f, (\lambda-\s)\f\> = \|f\|^2_{1,\lambda,\S}$, and
$\|\f\|_{-1,\lambda,\s}^2 = \sup_{\g {\rm \ local}}\{2\<\f,\g\> -
\|\g\|_{1,\lambda,\s}^2\}= \|f\|^2_{-1,\lambda,\S}$.

Also, in the following, it will be convenient to denote, when $B$
and its coefficient $\B$ are symmetric exclusion-type operators,
that $\|f\|^2_{1,\lambda, B}=\<f,(\lambda
 -B)f\>_\rho=\<\f,(\lambda - \B)\f\>=\|\f\|_{1,\lambda,\B}^2$ and
 $\|f\|_{-1,\lambda, B}= \sup_{g \ {\rm local}}\{2\<f,g\>_\rho -
 \|g\|^2_{1,\lambda, B}\} = \sup_{\g \ {\rm local}}\{2\<\f,\g\> -
 \|\g\|^2_{1,\lambda, \B}\} = \|\f\|^2_{-1,\lambda,\B}$.


\subsection{Some Variance Bounds and Comparisons} For a real local mean-zero
function $f$, $E_\rho[f]=0$, denote the variance
$$\sigma^2_t(f) \ = \ E_\rho\bigg[ \bigg(\int_0^t f(\zeta(s))ds\bigg)^2\bigg].$$
A well known upperbound on $\sigma^2_t(f)$, which connects with
$H_{-1}$ norms, and proved say in Proposition 6.1, appendix 1
\cite{KL}, is given in the next statement.
\begin{proposition}
\label{ub1} There is a universal constant $C_1$ such that for $t\geq
0$,
\begin{eqnarray*}
\sigma^2_t(f) &\leq& C_1t\<f, (t^{-1} - \L)^{-1} f\>_\rho \ .
\end{eqnarray*}
\end{proposition}


We now compare $\<f,(\lambda -\L)^{-1}f\>_\rho$ with other quadratic
forms depending on the dimension $d$.  Let $\L_{\rm nn}$ be the
reference process generator corresponding to nearest-neighbor jump
rates $p_{\rm nn}$ supported on standard vectors $\{\pm e_l\}$ of
$\Z^d$ where
$$
p_{\rm nn}(\pm e_l) \ = \ \left\{\begin{array}{rl} \max[\pm e_l
\cdot \sum jp(j), 0] & \
{\rm when \ }\pm e_l\cdot \sum jp(j) \neq 0\\
1& \ {\rm when\ } \pm e_l \cdot \sum jp(j) = 0\end{array}\right. $$
for $1\leq l\leq d$, and $p_{\rm nn}(z) = 0$ for $|z|\neq 1$. Note
that $s_{\rm nn}(z)=(p_{\rm nn}(z) +p_{\rm nn}(-z))/2>0$ for
$|z|=1$.

When $d=1$, define also operator $\N$ on local functions $f$ by
\begin{equation}\label{N}(\N f)(\zeta) \ = \
f(\zeta^{-1,1})-\f(\zeta),\end{equation} that is, the symmetric
exchange operator on bond connecting $-1$ and $1$.  Its coefficient
operator $\n$ defined on local functions $\f$ is then $(\n\f)(B) =
\f(B_{-1,1})-\f(B)$.

The next proposition, which indicates the $H_{-1}$ norm with respect
to $\L$ is on the same order as that for a nearest-neighbor dynamics
with the same drift, is Theorem 2.1 \cite{Scomp} for $d\geq 2$ and
proved by the proof of Theorem 2.2 \cite{Scomp} for $d=1$ (cf. Lemma
3.5 and p. 50 \cite{Scomp}).
\begin{proposition}
\label{comp} We have a constant $C=C(d,p)$, such that for $\lambda>0$
and local $f$
in $d\geq 2$,
$$C^{-1}\<f,(\lambda -\L_{\rm nn})^{-1}f\>_\rho \ \leq \ \<f,(\lambda
-\L)^{-1}f\>_\rho \ \leq \ C\<f,(\lambda -\L_{\rm nn})^{-1}f\>_\rho,
$$
and in $d=1$,
$$C^{-1}\<f,(\lambda -\L_{\rm nn} -\N)^{-1}f\>_\rho \ \leq \ \<f,(\lambda
-\L)^{-1}f\>_\rho \ \leq \ C\<f,(\lambda -\L_{\rm
nn}-\N)^{-1}f\>_\rho.
$$
\end{proposition}


 Let $\S_{\rm nn}$ and $\A_{\rm nn}$ be the
symmetric and anti-symmetric parts of $\L_{\rm nn}=\S_{\rm nn}
+\A_{\rm nn}$. Let also $\S_{\rm nn}^e$ and $\S_{\rm nn}^t$ be the
``environment'' and ``tagged-shift'' parts of $\S_{\rm nn} = \S_{\rm
nn}^e +\S_{\rm nn}^t$.  We denote also by $\a_{\rm nn}$ and $\s_{\rm
nn}^e$ the respective coefficients of $\A_{\rm nn}$ and $\S_{\rm
nn}^e$.


Recall the $H_1$ and $H_{-1}$ norm expressions
$\|\cdot\|_{\pm,\lambda,\B}$ for symmetric operators $\B$ at the end
of subsection 3.2.  The following bound allows us to bound $H_1$ and
$H_{-1}$ norms of the non-local ``tagged-shift'' operator $\S_{\rm
nn}^t$ in terms of the more manageable ``environment'' operator
$\S_{\rm nn}^e$. The proof is postponed to the last subsection of
this section.
\begin{proposition}
\label{tag_bounded_by_env} We have a constant $C=C(n,p)$ such that
for $\lambda>0$ and local $f$ with degree $n$ in $d\geq 2$,
$$\|f\|_{1,\lambda,\S_{\rm nn}^e} \ \leq \ \|f\|_{1,\lambda,\S_{\rm nn}} \ \leq \ C \|f\|_{\S_{\rm nn}^e},$$
and so consequently,
$$C^{-1} \|f\|_{-1,\lambda,\S_{\rm nn}^e} \ \leq \
\|f\|_{-1,\lambda,\S_{\rm nn}} \ \leq \|f\|_{-1,\lambda,\S_{\rm
    nn}^e}.$$
In $d=1$, the inequalities hold with $\S_{\rm nn}^e$ replaced by
$\S^e_{\rm nn} +\N$.
\end{proposition}
%


\subsection{``Extended'' Coefficient Functions} To aid later computations,
we now extend the underlying space $\Z^d\setminus\{0\}$ to $\Z^d$.
We concentrate on dimension $d\leq 2$ for simplicity. Let
$\bar{\Epsilon}_d$ be the set of finite subsets of $\Z^d$, and let
$\bar{\Epsilon}_{d,n}$ be those subsets of $\Z^d$ with cardinality
$n$.
 Let also $\bar{\C}_{d,n}$ denote the collection of functions on $\bar{\Epsilon}_{d,n}$.

 For $n\leq 2$, let $\f\in \C_{d,n}$
be a coefficient function. We now give extensions $\f_{\rm ext}$ and
$\f_\odot$ belonging to $\bar{\C}_{d,n}$; we also give an
``inverse'' of the $\odot$ extension, namely $\g_{\rm res}$, which
restricts $\g\in \bar{\C}_{d,n}$ to $\C_{d,n}$. In addition, we
define some related operators, an innerproduct, and norms, acting on
these functions. \vskip .2cm

{\bf Extension $\f_{\rm ext}$.} This extension assigns to sets $B\ni
0$ the ``local'' average of ``nearest-neighbor'' sets and is well
suited for later comparisons of Dirichlet forms over $\Z^d\setminus
\{0\}$ and $\Z^d$ (cf. Proposition \ref{odot_bounds}).  More
precisely, when $n=1$, let
$$
\f_{{\rm ext}}(\{x\}) \ =\ \left\{\begin{array}{rl} \f(\{x\}) & \ {\rm  for\ }x\in \Z^d\setminus\{0\}\\
\frac{1}{2d}\sum_{|z|=1} \f(\{z\})& \ {\rm for \
}x=0.\end{array}\right. $$
 When $n=2$, for distinct $x,y\in
\Z^d\setminus\{0\}$, let $\f_{{\rm ext}}(\{x,y\}) =\f(\{x,y\})$, and
$$\f_{{\rm ext}}(\{0,y\})  \ = \ \left\{\begin{array}{rl}
\frac{1}{2d-1} \sum_{z\neq y\atop |z|=1}f(\{z,y\}) &\ {\rm when \ }
|y|=1\\
\frac{1}{2d} \sum_{|z|=1} \f(\{z,y\})& \ {\rm when \ }|y|\geq 2.

\end{array}\right.
$$

\vskip .1cm

{\bf Extension $\f_\odot$.} This type of extension vanishes on sets
involving the origin and allows $H_{-1}$ norm comparisons over
$\Z^d\setminus\{0\}$ and $\Z^d$ (cf. Proposition \ref{odot_bounds}).
Let
$$
\f_\odot(B) =
\left\{\begin{array}{rl} \f(B)&
    \ {\rm when \ }B\in \Epsilon_{d,n}\\
0& \ {\rm otherwise.}\end{array} \right. $$ \vskip .1cm

{\bf Restriction $\g_{\rm res}$.}  For $\g \in \bar{\C}_{d,n}$, let
$\g_{\rm res}\in \C_{d,n}$ be the restriction of $\g$ to subsets
$B\in \Epsilon_{d,n}$.  This restriction is useful in extending
operators with respect to $\Z^d\setminus\{0\}$ to underlying space
$\Z^d$ (cf. definition of $\bar{\a}_{{\rm nn}; n,m}$ below).


\vskip .2cm

{\bf Operator $\s_{\rm ext}$.}  Recall operator $\S_{\rm nn}^e$ and
its coefficient form $\s_{\rm nn}^e$ from subsection 3.3.  We now
extend $\s^e_{\rm nn}$ on local $\C_{d,n}$ functions to $\s_{\rm
ext}$ acting on local $\bar{\C}_{d,n}$ functions in the usual way,
namely transitions are now allowed into the origin. Define the
nearest-neighbor operator, acting on local $\g\in\bar{\C}_{d,n}$, by
$$
(\s_{\rm ext}\g)(B) \ = \ \sum_{|i-j|=1\atop i,j\in \Z^d}
\bigg(\g(B_{i,j}) - \g(B) \bigg).$$

\vskip .1cm

{\bf Operator $\bar{\a}_{{\rm nn}; n,m}$.}  Recall operator $\A_{\rm
nn}$ and its coefficient form $\a_{\rm nn}$ in subsection 3.3. With
respect to $\a_{{\rm nn}; n,m}$, the part of $\a_{\rm nn}$ which
takes degree $n$ functions to degree $m$, define on local
$\g\in\bar{\C}_{d,n}$ that
$$
(\bar{\a}_{{\rm nn}; n,m}\g)(B) \ = \
\left\{\begin{array}{rl}
(\a_{{\rm nn}; n,m}{\g_{\rm res}})(B) & \ {\rm when \ } B\subset \Z^d\setminus\{0\}\\
0& \ {\rm otherwise.}\end{array}\right . $$ \vskip .1cm

{\bf Extended Innerproduct and Norms.} The innerproduct naturally
extends to $L^2$ functions in $\bar{\C}_{d,n}$:
$$\<\f,\g\>_{\rm ext} \ = \ \sum_{|B|=n\atop B\subset \Z^d}\f(B)\g(B).$$
Also, $H_1$ and $H_{-1}$ norms of $\f\in \bar{\C}_{d,n}$, with
respect to $\s_{\rm ext}$, are defined for $\lambda>0$:
\begin{eqnarray}  \label{extform}
\|\f\|^2_{1,\lambda, \s_{\rm ext}} &=& \<\f,(\lambda -\s_{\rm
ext})\f\>_{\rm ext}\nonumber\\
&=& \lambda \<\f,\f\>_{\rm ext} + \frac{1}{2}\sum_{|B|=n\atop
B\subset \Z^d}\sum_{|i-j|=1\atop i,j\in
\Z^d}(\f(B_{i,j})-\f(B))^2\\
 \|\f\|^2_{-1, \lambda, \s_{\rm ext}} &
= & \sup_{\g\ {\rm local \ on \ }\bar{\Epsilon}_d}
\big\{2\<\f,\g\>_{\rm ext} - \<\g, (\lambda - \s_{\rm ext})
\g\>_{\rm
    ext}\big\}.\nonumber
    \end{eqnarray}

In addition, we have the following useful bounds which relate
further the various extensions.
\begin{lemma}
\label{g'bounds} For $0<\lambda\leq 1$, we have a constant $\C(d)$
such that for $\g\in \C_{d, 1}$ and any extension $\g'\in
\bar{\C}_{d,1}$, $$
\|\g_{\rm ext}\|^2_{1,\lambda,\s_{\rm ext}} \ \leq \
C\bigg[\|\g'\|^2_{1,\lambda,\s_{\rm ext}} + |\g'(\{0\}) -
\sum_{|z|=1}\g'(\{z\})|^2\bigg].$$
\end{lemma}

{\it Proof.} Note first $\g_{\rm ext} = \g' +[\sum_{|z|=1}\g'(\{z\})
- \g'(\{0\})]\omega_0$ where $\omega_0\in \bar{\C}_{d,1}$ and
$\omega_0(\{x\}) = 1$ for $x=0$ and vanishes otherwise.  Then,
$$\|\g_{\rm ext}\|^2_{1,\lambda,\s_{\rm ext}} \leq 2\|\g'\|^2_{1,\lambda,\s_{\rm ext}}
+ 2\big[\sum_{|z|=1}\g'(\{z\}) -
\g'(\{0\})\big]^2\|\omega_0\|^2_{1,\lambda,\s_{\rm ext}}.$$ By
calculation, using (\ref{extform}),
$\|\omega_0\|^2_{1,\lambda,\s_{\rm ext}} \leq \lambda + C$ and so
the result follows. \qed

  Recall symmetric operators
$\S^e_{\rm nn}$ and $\N$, and their coefficients $\s^e_{\rm nn}$ and
$\n$ from subsection 3.3, and $H_1$ and $H_{-1}$ norm expressions
 $\|\cdot\|_{\pm, \lambda, \B}$ for symmetric operators $\B$ at the
end of subsection 3.2.
\begin{proposition}
\label{odot_bounds}  For $n\leq 2$ and $\lambda>0$, we have a
constant $C=C(d,n,p)$ such that for $\f\in {\C}_{d,n}$ in $d=2$,
\begin{equation}
\label{lb}C^{-1}\|\f\|_{1,\lambda,\s_{\rm nn}^e}\ \leq \ \|\f_{{\rm
ext}}\|_{1,\lambda, \s_{\rm
    ext}} \ \leq \ C \|\f\|_{1,\lambda, \s_{\rm
    nn}^e}\end{equation}
and
$$\|\f\|_{-1,\lambda, \s_{\rm nn}^e} \ \leq \ C\|\f_\odot\|_{-1,\lambda, \s_{\rm ext}}\ .$$
In $d=1$, the inequalities hold with operator $\s^e_{\rm nn}$
replaced by $\s^e_{\rm nn} + \n$.
\end{proposition}
We postpone the proof to the last subsection of this section.

\subsection{``Free Particle'' Bounds.} For later detailed analysis, it will be helpful
to ``remove the hard-core exclusion.''  In other words, we want to
get equivalent bounds in terms of operators which govern completely
independent or ``free'' motions.  We follow the treatment of
Bernardin \cite{B} with respect to occupation times.

\vskip .2cm

{\bf ``Free Particle'' Generator $\s_{\rm free}$.} Let
$\upsilon_{d,n} = (\Z^d)^n$ and consider $n$ independent random
walks with symmetric nearest-neighbor symmetric jump rates on $\Z^d$
for $d\geq 1$. The process $x_t = (x^1_t,\ldots,x^n_t)$ evolves on
$\upsilon_{d,n}$ and has generator $\s_{\rm free}$ acting on local,
namely finitely supported, functions on $\upsilon_{n,d}$,
$$ (\s_{\rm free}\f)(x) \ = \
\frac{1}{2d}\sum_{1\leq j\leq n \atop |z| = 1}
\bigg(\phi(x+z\omega_j)-\phi(x)\bigg)$$
  where $z\omega_j = (0,\ldots,0,z,0,\ldots,0)$ is the state with
  $z$ in the $j$th place.

\vskip .2cm {\bf Free Innerproduct and Norms.} With respect to local
functions on $\upsilon_{d,n}$, define $$\<\phi,\psi\>_{{\rm free}} =
\frac{1}{n!}\sum_{x\in \upsilon_{d,n}} \phi(x)\psi(x).$$

Define also, for $\lambda>0$, $H_{1,\lambda}$ and $H_{-1,\lambda}$ norms
$\|\phi\|^2_{1,\lambda,{\rm free}} = \<\phi, (\lambda -\s_{\rm
free})\phi\>_{{\rm free}}$ and
$$\|\phi\|^2_{-1,\lambda, {\rm
free}} \ = \ \sup_{\psi \ {\rm local \ on \
}\upsilon_{d,n}}\big\{2\<\phi,\psi\>_{{\rm free}}-
\|\psi\|^2_{1,\lambda,{\rm free}}\big\}.$$
\vskip .1cm

{\bf Extension $\f_{\rm free}$.} Let $\G_n\subset \upsilon_{d,n}$ be
those points whose coordinates are distinct.  For a function $\f\in
\bar{\C}_{d,n}$, define the natural extension to $\upsilon_{d,n}$ by
$$
\f_{{\rm free}}(x) \ =\ \f(U)$$ where $U$ is the set formed from
coordinates of $x\in \upsilon_{d,n}$. Note $\f_{{\rm free}}$ is
supported on $\G_n$.

\vskip .2cm

{\bf Extension $\tilde{\f}$.} We now give an extension $\tilde{\f}$
on $\upsilon_{d,n}$ which allows some $H_1$ and $H_{-1}$ norm
comparsions (cf. Proposition \ref{free}).  Let $\tau$ be the arrival
time into $\G_n$,
$$\tau \ = \ \inf\big\{t\geq 0: x_t \in \G_n\big\}.$$
Then, for $\f\in \bar{\C}_{d,n}$, define for $x\in \upsilon_{n,d}$
that
$$
\tilde{\f}(x) \ = \ E_x[\f_{\rm free}(x_\tau)].$$

%
%
\vskip .2cm {\bf Free Bounds and Relations.}
 The next result relates $\s_{\rm ext}$ and $\s_{\rm free}$ with respect to
 $H_{1}$ and $H_{-1}$ norms of $\f$ and $\tilde{\f}$,
 and is a part of Theorems 3.1 and 3.2
\cite{B}.
\begin{proposition}
\label{free} We have, for a constant $C=C(d,n,p)$, $\lambda>0$, and $\f\in
\bar{\C}_{d,n}$, that
$$C^{-1}\|\tilde{\f}\|_{1,\lambda,{\rm free}} \ \leq \
\|\f\|_{1,\lambda,\s_{\rm ext}}
\ \leq C\|\tilde{\f}\|_{1,\lambda, {\rm free}}.$$ Also,
$$ \|\f\|_{-1,\lambda, \s_{\rm ext}} \ \leq \ C\| 1_{\G_n}\tilde{\f}\|_{-1, \lambda,
  {\rm free}}.$$
\end{proposition}

The following relations, which follow from straightforward
manipulations, will also be useful.
\begin{lemma}\label{free_relation}
Let $\g\in \C_{d,1}$ be a local function, and let $\g'\in
\bar{\C}_{d,1}$ be any extension. Then, for $x\in \Z^d$
$$\widetilde{\g'}(x) \ = \ \g'_{\rm free}(x), \ \ \ {\rm and \
}1_{\G_1}\widetilde{(\a_{{\rm nn;}1,1}\g)_\odot}(x) \ =\
(\bar{\a}_{{\rm nn;}1,1}\g')_{\rm free}(x).$$ Also, for $x,y\in \Z^d$,
$$ 1_{\G_2}\widetilde{(\a_{{\rm nn;}1,2}\g)_\odot}(x,y) \ = \
(\bar{\a}_{{\rm nn;}1,2}\g')_{\rm
  free}(x,y).$$
  \end{lemma}
 \vskip .1cm

{\bf Fourier Transform Expressions.} It will be convenient to
express ``free'' $H_{1,\lambda}$ and $H_{-1,\lambda}$ norms in terms
of Fourier transforms.  Let $\psi$ be a local function on
$\upsilon_{d,n}$ and let $\widehat{\psi}$ be its Fourier transform
$$\widehat{\psi}(s_1,\ldots,s_n) = \frac{1}{\sqrt{n!}} \sum_{x\in \upsilon_{d,n}}
e^{2\pi i (x_1\cdot s_1 + \cdots + x_n\cdot s_n)}\psi(x)$$ where
$s_1,\ldots,s_n\in [0,1]^d$. Compute
$$\widehat{\s_{\rm free}\psi}(s_1,\ldots,s_n) \ = \ -\bigg[\sum_{j=1}^n
\theta_d(s_j)\bigg] \widehat{\psi}(s_1,\ldots,s_n)$$ where
$\theta_d(u) = (2/2d)\sum_{z\in \Z^d\atop |z|=1} \sin^2(\pi(u\cdot
z))=(2/d)\sum_{j=1}^n \sin^2(\pi u_j)$. Hence, we have
$$\|\psi\|^2_{1,\lambda, {\rm free}} \ = \ \int_{{s\in ([0,1]^d)^n\atop s=(s_1,\ldots,s_n)}}
\bigg(\lambda + \sum_{j=1}^n
\theta_d(s_j)\bigg)|\widehat{\psi}(s_1,\ldots,s_n)|^2 ds$$ and
$$\|\psi\|^2_{-1,\lambda, {\rm free}} \ = \ \int_{{s\in ([0,1]^d)^n\atop s=(s_1,\ldots,s_n)}}
\frac{|\widehat{\psi}(s_1,\ldots,s_n)|^2}{\lambda + \sum_{j=1}^n
\theta_d(s_j)} ds.$$

\subsection{Putting Bounds Together} We now incorporate the
previous bounds into a single statement.

\begin{proposition}
\label{mainprop} In $d\leq 2$, for local degree one functions $f\in
\C_{d,1}$, we have a constant $C=C(d,p,\rho)$ such that for $t\geq
1$,
\begin{eqnarray*} \sigma^2_t(f)/t & \leq & C\inf_{\g\in \bar{\C}_{d,1}, \atop \ \ {\rm
local}}\bigg\{\|(\f_\odot)_{\rm free} - (\bar{\a}_{{\rm
nn};1,1}\g)_{\rm free}\|^2_{-1,t^{-1},{\rm
free}}+ \|\g_{\rm free}\|^2_{1,t^{-1},{\rm free} }\\
&&\ \ \ \ \ \ \ \ \ \ \ \ \ \ \ + \|(\bar{\a}_{{\rm nn};1,2}\g)_{\rm free}\|^2_{-1,t^{-1},{\rm
free}}  +|\g_{\rm
free}(0)-\sum_{|z|=1}\g_{\rm free}(z)|^2\bigg\}.
\end{eqnarray*}
\end{proposition}

{\it Proof.} In the following, the constant $C=C(d,p,\rho)$ can
change from line to line. We have, in sequence, from Propositions
\ref{ub1}, \ref{comp}, \ref{infvarexpression}, \ref{tag_bounded_by_env} and
\ref{odot_bounds}, when $d=2$ that
\begin{eqnarray*}
\sigma^2_t(f)/t &\leq& C\<f,(t^{-1} -\L)^{-1}f\>_\rho\\
&\leq& C\<f,(t^{-1} -\L_{\rm nn})^{-1}f\>_\rho\\
&=& C\inf_{g\ {\rm local}}\{\|f-\A_{\rm nn}g\|^2_{-1,t^{-1},\S_{\rm
nn}} +
\|g\|^2_{1,t^{-1},\S_{\rm nn}}\bigg\}\\
 &\leq& C\inf_{g\ {\rm local}}\bigg\{\|f-\A_{\rm nn}g\|^2_{-1,t^{-1},\S_{\rm nn}^e}
 + \|g\|^2_{1,t^{-1},\S_{\rm nn}^e}\bigg\}\\
&\leq& C\inf_{\g\ {\rm local}}\bigg\{\|\f_\odot - (\a_{\rm
nn}\g)_\odot\|^2_{-1,t^{-1},\s_{\rm ext}} + \|\g_{\rm
ext}\|^2_{1,t^{-1},
  \s_{\rm ext}}\bigg\}.
\end{eqnarray*}
When $d=1$, in the fourth line of the sequence above, $\S_{\rm nn}^e$ is replaced by
$\S_{\rm nn}^e +\N$.


The last
infimum, by first restricting to $\g\in \C_{d,1}$ and Schwarz inequality, second using
Lemma \ref{g'bounds} to estimate $\|\g_{\rm ext}\|_{1,t^{-1},\s_{\rm
ext}}^2$ in terms of $\|\g'\|_{1,t^{-1},\s_{\rm ext}}^2$ for $\g'\in
\bar{\C}_{d,1}$, and then third applying Proposition \ref{free} and
Lemma \ref{free_relation} to estimate in terms of ``free'' norms on
local functions in $\bar{\C}_{d,1}$, is further bounded by twice
\begin{eqnarray*}
&&\inf_{\g \in \C_{d,1},\atop \ \  \ {\rm
local}}\bigg\{\|\f_\odot - (\a_{{\rm nn};1,1}\g)_\odot\|^2_{-1,t^{-1},\s_{\rm
ext}}\nonumber\\
&&\ \ \ \ \ \ \ \ \ \ \ \ \ \ \ \ \ \ \ \ \  + \|(\a_{{\rm nn};1,2}\g)_\odot\|^2_{-1,t^{-1},\s_{\rm ext}} +
\|\g_{\rm ext}\|^2_{1,t^{-1},
  \s_{\rm ext}}\bigg\}\nonumber\\
&&\ \  \leq C\inf_{\g'\in \bar{\C}_{d,1}, \atop \ \ \ {\rm
local}}\bigg\{\|(\f_\odot)_{\rm free} - (\bar{\a}_{{\rm nn};1,1}\g')_{\rm
free}\|^2_{-1,t^{-1},{\rm
free}} \\
&&\ \ \ \ \ \ \ \ \ \ \ \ \ \ \ \ \ \  + \|(\bar{\a}_{{\rm nn};1,2}\g')_{\rm
free}\|^2_{-1,t^{-1},{\rm free}}\nonumber\\
&&\ \ \ \ \ \ \ \ \ \ \ \ \ \ \ \ \ \ + \|\g'_{\rm
free}\|^2_{1,t^{-1},{\rm free}} + |\g'_{\rm
free}(0)-\sum_{|z|=1}\g'_{\rm free}(z)|^2\bigg\}.\nonumber
\end{eqnarray*}
\qed

\subsection{Proofs of Propositions \ref{tag_bounded_by_env} and
\ref{odot_bounds}}
\hfill

 {\it Proof of Proposition
\ref{tag_bounded_by_env}.} 
The $H_1$ lower bound follows as
$$\<f,(-\S_{\rm nn})f\>_\rho \ =\
\<f,(-\S_{\rm nn}^e)f\>_\rho +\<f,(-\S_{\rm nn}^t)f\>_\rho$$
 and
$$\<f,(-\S_{\rm nn}^t)f\>_\rho \ =\
\frac{1}{2}\sum_{|z|=1}s_{\rm
nn}(z)E_\rho[(1-\zeta_z)(f(\tau_z\zeta)-f(\zeta))^2]\ \geq \ 0.$$

For the $H_1$ upper bound, note
\begin{eqnarray*}
E_\rho[(1-\zeta_z)(f(\tau_z\zeta)-f(\zeta))^2] & \leq &
E_\rho[(f(\tau_z \zeta)-f(\zeta))^2] \\
& =& \sum_{B\subset
\Z^d\setminus\{0\}\atop |B|=n} (\f(\tau_{-z}B) -
\f(B))^2,\end{eqnarray*} and by the proof of Lemma 5.1 \cite{lan-olla-var},
\begin{eqnarray}
 \sum_{B\subset\Z^d\setminus\{0\}\atop |B|=n}
(\f(\tau_{-z}B) - \f(B))^2& \leq & C_z
n\sum_{B\subset\Z^d\setminus\{0\}\atop |B|=n}\sum_{i\sim
j}\ (\f(B_{i,j}) - \f(B))^2 \label{rightsideexp}\end{eqnarray} where
$C_z$ is a constant depending on $z$,
and $i\sim j$ means a ``neighboring'' pair $i,j\in
\Z^d\setminus\{0\}$ with $|i-j|=1$, or also $(i,j)=(1,-1)$ and
$(-1,1)$ when $d=1$.  Also
$$ \sum_{B\subset\Z^d\setminus\{0\}\atop |B|=n} \sum_{i\sim  j} \ (\f(B_{i,j})-\f(B))^2
\ \leq \  \left\{\begin{array}{rl}
C'\<f,(-\S^e_{\rm nn})f\>_\rho & \ {\rm when \ }d\geq 2\\
C'\<f,(-\S^e_{\rm nn}-\N)f)\>_\rho & \ {\rm when \ }d=1\end{array}\right. 
$$
where $C' = C'(s_{\rm nn})$.
The $H_1$ estimates in the proposition follow now by adding
over $|z|=1$. Also, the $H_{-1}$ bounds are deduced from the $H_1$
bounds through simple estimates with the definition of
$\|f\|^2_{-1,\lambda, \B}$ (cf. subsection 3.2). \qed
\vskip .1cm

{\it Proof of Proposition \ref{odot_bounds}.} We prove the statement
for $d=2$, and mention at the end modifications for $d=1$. In the
following, $C=C(n,p)$ denotes a constant which can change from line
to line. The lowerbound inequality in (\ref{lb}) follows from
overcounting:
\begin{eqnarray*} \<\f,(-\s_{\rm
nn}^e)\f\>&=& \frac{1}{2}\sum_{B\subset\Z^2\setminus\{0\}\atop
|B|=n}\sum_{|i-j|=1\atop i,j\neq
  0}(\f(B_{i,j})-\f(B))^2s_{\rm nn}(j-i)\\
  &\leq& \frac{C}{2}\sum_{B\subset \Z^2\atop
|B|=n}\sum_{|i-j|=1\atop i,j \in\Z^2} (\f_{\rm ext}(B_{i,j}) -
\f_{\rm
  ext}(B))^2 \ = \ C\<f_{{\rm ext}}, (-\s_{\rm ext}) \f_{\rm
ext}\>_{\rm ext}; \end{eqnarray*} also, we have $\|\f\|_0^2 \leq
\<\f_{\rm ext},\f_{\rm ext}\>_{\rm ext}$.

For the upperbound in (\ref{lb}), as $s_{\rm nn}(z)>0$ for $|z|=1$,
we have
\begin{eqnarray*}
\<\f,(-\s_{\rm nn}^e)\f\>
  & \geq & C\sum_{B\subset\Z^2\setminus\{0\}\atop
|B|=n}\sum_{|i-j|=1\atop i,j\in
\Z^2\setminus\{0\}}(\f(B_{i,j})-\f(B))^2\end{eqnarray*} 
and so 
\begin{eqnarray}
\label{secondterm} \<f_{{\rm ext}}, (-\s_{\rm ext}) \f_{\rm
ext}\>_{\rm ext} &\leq& C\<\f,(-\s_{\rm nn}^e)\f\>\nonumber\\
&&\  + C\sum_{|i-j|=1\atop i,j\in \Z^2}\sum_{B \ {\rm or \
}B_{i,j}\ni 0 \atop |B|=n}(\f_{\rm ext}(B_{i,j}) - \f_{\rm
ext}(B))^2.
\end{eqnarray}
When $n=1$, the last term of (\ref{secondterm}) is on order
$$\sum_{|z|=1}(\f_{\rm
ext}(\{0\})-\f_{\rm ext}(\{z\}))^2 \ = \
\frac{1}{16}\sum_{|w|,|z|=1}(\f(\{w\})-\f(\{z\}))^2 \ \leq \
C\<\f,(-\s_{\rm nn}^e)\f\>.$$ Here, for the last inequality, we
build a path from $w_0=e_1$ to $w_1=e_1+e_2$ to $w_2=e_2$ and so on
to $w_7=e_1-e_2$ back to $w_8=e_1$, and bound each of the finite
number of terms $(\f(\{w\})-\f(\{z\}))^2 \leq 8\sum_{i=0}^7
(\f(\{w_i\})-\f(\{w_{i+1}\}))^2 \leq C\<\f,(-\s_{\rm nn}^e)\f\>$.

When $n=2$, the last sum in (\ref{secondterm}) is on order
\begin{eqnarray*}&&\sum_{y\neq 0}\bigg[\sum_{|z|=1\atop y+z\neq 0}(\f_{\rm
ext}(\{0,y+z\})-\f_{\rm ext}(\{0,y\}))^2 + \sum_{|z|=1\atop
z\neq y}(\f_{\rm ext}(\{z,y\}) - \f_{\rm ext}(\{0,y\}))^2\bigg]\\
&&\  = \sum_{|y|\geq 2}\sum_{|z|=1}\bigg[(\f_{\rm
ext}(\{0,y+z\})-\f_{\rm
ext}(\{0,y\}))^2 +(\f_{\rm ext}(\{z,y\}) - \f_{\rm ext}(\{0,y\}))^2\bigg] \\
&&\ \ \ \ \ \ \ \ \ \ \ \ \ \ + \ {\rm finite \ number \ of \
remaining \ terms}.
\end{eqnarray*}
The first line is straightforwardly bounded by $C\<\f,(-\s_{\rm
nn}^e)\f\>$. The remaining finite number of terms are handled as
follows:  For $|y|=|z|=1$, the terms with $y+z\neq 0$ are bounded
\begin{eqnarray*}\bigg(\f_{\rm ext}(\{0,y+z\}) - \f_{\rm
ext}(\{0,y\})\bigg)^2 & =&
\bigg(\frac{1}{4}\sum_{|x|=1}\f(\{x,y+z\})
-\frac{1}{3}\sum_{|w|=1\atop w\neq y} \f(\{w,y\})\bigg)^2 \\
& \leq & C\<\f,(-\s_{\rm nn}^e)\f\>\end{eqnarray*} and the terms
with $|y|=|z|=1$ and $z\neq y$ are bounded
$$\bigg(\f_{\rm ext}(\{z,y\}) - \f_{\rm ext}(\{0,y\})\bigg)^2
\ = \  \bigg(\f(\{z,y\}) - \frac{1}{3}\sum_{|w|=1 \atop w\neq
y}\f(\{w,y\})\bigg)^2 \ \leq \ C\<\f,(-\s_{\rm nn}^e)\f\>$$ through
similar arguments using the path built in the $n=1$ case.

Also, more directly, $\<\f_{\rm
  ext},\f_{\rm ext}\>_{\rm ext} \leq C\|\f\|_0^2$ to finish the
  upperbounds in the first statement of the proposition.

For the second statement after (\ref{lb}), write
\begin{eqnarray*}
\|\f\|_{-1,\lambda, \s_{\rm nn}^e}^2 &=&\sup_{\phi \ {\rm
local}}\big\{2\<\f,\phi\> -
\|\phi\|^2_{1,\lambda,\s_{\rm nn}^e}\big\}\\
&=&\sup_{\phi\in \C_{2,n}\ {\rm local}}\big\{2\<\f,\phi\> -
\|\phi\|^2_{1,\lambda,\s_{\rm nn}^e}\big\}. \end{eqnarray*} The last
step follows as for $\f\in \C_{2,n}$ with $\phi = \sum_m\phi_m$
decomposed in degrees, $\<\f,\phi\> = \<\f, \phi_n\>$ and as
$\s_{\rm nn}^e$ preserves degrees, $\|\phi\|^2_{1,\lambda,\s_{\rm
nn}^e} = \sum_{m}\|\phi_m\|^2_{1,\lambda,\s_{\rm nn}^e}$; so one
does best by choosing $\phi = \phi_n$.

Continuing, as $\<\f,\phi\> = \<\f_\odot, \phi_{\rm ext}\>_{\rm
ext}$ and using the proved lowerbound in (\ref{lb}),
$\|\f\|^2_{-1,\lambda,\s^e_{\rm nn}}$ is bounded above by
$$\sup_{\phi\in \C_{2,n} \ {\rm local}}\big\{2\<\f_\odot, \phi_{{\rm ext}}\>_{\rm ext} - C^{-1}\|\phi_{{\rm
    ext}}\|^2_{1,\lambda,\s_{\rm ext}}\big\}\ \leq \ C\|\f_\odot\|^2_{-1,\lambda, \s_{\rm
ext}}.
$$

The modifications for $d=1$ take advantage of inequalities
\begin{eqnarray*} \<\f,(-\n)\f\> & = &
\frac{1}{2}\sum_{B\subset \Z\setminus\{0\}\atop
|B|=n}(\f(B_{1,-1})-\f(B))^2\\
&= &\frac{1}{2}\sum_{B\subset \Z\setminus\{0\}\atop
|B|=n}(\f_{\rm ext}(B_{1,-1})-\f_{\rm ext}(B))^2\\
&\leq & \frac{1}{2}\sum_{B\subset \Z\atop |B|=n}(\f_{\rm
ext}(B_{1,-1})-\f_{\rm ext}(B))^2\\
&\leq &C\sum_{B\subset \Z\atop |B|=n}\bigg[(\f_{\rm
ext}(B_{1,0})-\f_{\rm ext}(B))^2 +(\f_{\rm ext}(B_{0,-1})-\f_{\rm
ext}(B))^2\bigg]
\end{eqnarray*}
which hold as $B_{-1,1} = ((B_{1,0})_{0,-1})_{1,0}$ and by applying
Schwarz inequality.  The arguments are now similar to those in
$d=2$.
 \qed
\sect{Proof of Theorem \ref{maintheorem1}}
First, by (\ref{martdecomp}) and that quadratic variation
$E_\rho[|M(t)|^2] = (1-\rho)t\sum_j |j|^2p(j) =O(t)$, we need only
bound
$$E_\rho[|A(t)|^2] \ \leq \ \sum |j|^2 \sigma^2_t(\rho-\zeta_j)p(j) \ = \
O(t).$$ Clearly, it is sufficient to show that $\sigma^2_t(\rho -
\zeta_{j_0}) = O(t)$ for $j_0\in \Z^2\setminus\{0\}$ and $t\geq 1$.




To accomplish this, through Proposition \ref{mainprop}, it will be
useful to compute, for a local function $\g\in \bar{\C}_{d,1}$,
Fourier transforms $\widehat{(\bar{\a}_{{\rm nn};1,1}\g)_{\rm
free}}$ and $\widehat{(\bar{\a}_{{\rm nn};1,2}\g)_{\rm free}}$ where
$\bar{\a}_{{\rm nn};n,m}$ are the nearest-neighbor operators defined
in subsection 3.4. When $d=2$, let $a_1=a_{\rm nn}(e_1)$ and
$a_2=a_{\rm nn}(e_2)$, and when $d=1$ let $a_1 = a(1)$. Note, by the
assumption $\sum jp(j)\neq 0$, that $a_1^2 +a_2^2>0$ in $d=2$ and
$|a_1|>0$ in $d=1$.

Let $\gamma(r) = e^{2\pi i r}-e^{-2\pi i r} = 2i\sin(2\pi r)$ for
$r\in [0,1]$.
The following proposition is proved in section 5.

\begin{proposition}
\label{FT} In $d\leq 2$, for local $\g\in \bar{\C}_{d,1}$ and a
constant $C=C(d,p,\rho)$,
\begin{eqnarray*}
\widehat{(\bar{\a}_{{\rm nn};1,1}\g)_{\rm free}}(v) &=&
\rho\bigg[\sum_{i=1}^d a_i\gamma(v_i)\bigg ]\widehat{\g_{\rm
free}}(v) + \delta_0(v)
\end{eqnarray*}
where $|\delta_0(v)| \leq \kappa(v)\sum_{|z|\leq 1}|\g_{\rm
free}(z)|$ and $\kappa(v)$ is a bounded function such that
$$\kappa(v)^2 \ \leq \ C|v-z|^2$$
as $v\rightarrow z$ for $z=(0,0)$, $(0,1)$, $(1,0)$ and $(1,1)$ in
$d=2$, and $z=0$ and $1$ in $d=1$. Also,
\begin{eqnarray*}
\sqrt{2}\widehat{(\bar{\a}_{{\rm nn};1,2}\g)_{\rm free}}(v,w) &=&
2\beta_\rho\bigg[\sum_{i=1}^d a_i\gamma(v_i+w_i) + \alpha_d(v,w)\bigg]\widehat{\g_{\rm free}}(v+w)\\
&&\ \ + \beta_\rho
\bigg[-\sum_{i=1}^da_i\gamma(v_i)+\alpha_d(v,w)\bigg]\widehat{\g_{\rm
free}}(v)\\
&&\ \ +\beta_\rho\bigg[-\sum_{i=1}^da_i\gamma(w_i) +
\alpha_d(w,v)\bigg]\widehat{\g_{\rm free}}(w) +\delta_1(v,w)
\end{eqnarray*}
where, for $r,s\in [0,1]^d$, \begin{eqnarray*} \alpha_d(r,s) &=&
\sum_{i=1}^d a_i\bigg[\gamma(r_i)+\gamma(s_i) -
\gamma(r_i+s_i)\bigg]\end{eqnarray*}
%
and $|\delta_1(v,w)|\leq \kappa(v,w)\sum_{|z|\leq 1}|\g_{\rm
free}(z)|$ and $\kappa(v,w)$ is a bounded function such that
$$\kappa(v,w)^2 \ \leq \
C[|v-z_1|^2 + |w-z_2|^2]$$ as $(v,w)\rightarrow (z_1,z_2)$ for
$z_1,z_2 = (0,0)$,$(0,1)$, $(1,0)$ and $(0,1)$ in $d=2$, and
$z_1,z_2 = 0$ and $1$ in $d=1$.

%
\end{proposition}

Let now $f(\zeta) = \rho-\zeta_{j_0}$. As $$(\f_\odot)_{\rm
free}(z) \ =\ \left\{\begin{array}{rl}-\beta_\rho & z=j_0\\
0&\ {\rm otherwise},\end{array}\right.$$ we calculate
$$\widehat{(\f_\odot)_{\rm free}}(v) \ =\
-\beta_\rho e^{2\pi i(j_0\cdot v)} \ =\ -\beta_\rho +\delta_2(v)$$
where $\delta_2(v) = -\beta_\rho(e^{2\pi i(j_0\cdot v)}-1)$ and so
$|\delta_2(v)|^2\leq C|v-z|^2$ as $v\rightarrow z$ for $z=
(0,0),(0,1)$, $(1,0)$, and $(1,1)$ in $d=2$, and $z=0$ and $1$ in
$d=1$.

We now apply Propositions \ref{mainprop} and \ref{FT}.
Write, for local $\g \in \bar{\C}_{d,1}$ and $\lambda = t^{-1}$, in
Fourier expression (cf. subsection 3.5), that $\sigma^2_t(f)/t$ is
less than
\begin{eqnarray}
\label{lastterm} 
&&{2}\int_{[0,1]^d} \frac{|-\beta_\rho-\rho \big [ \sum_{i=1}^d
a_i\gamma(v_i)\big]\widehat{\g_{\rm free}}(v)|^2}{\lambda +
\theta_d(v)} + (\lambda +\theta_d(v))|\widehat{\g_{\rm free}}(v)|^2\
dv
\\
&&+ 2\int_{[0,1]^d} \frac{|\delta_0(v) +\delta_2(v)|^2}{\lambda
+\theta_d(v)}dv \ + \ \bigg|\g_{\rm free}(0)-\sum_{|z|=1}\g_{\rm
free}(z)\bigg|^2\nonumber \\
&& + \frac{3}{2}\beta_\rho^2\int_{([0,1]^d)^2} \frac{dvdw}{\lambda +
\theta_d(v) +\theta_d(w)}\nonumber \\
&&\ \  \times \bigg|2\widehat{\g_{\rm free}}(v+w)\sum_{i=1}^d
a_i\gamma(v_i+w_i) -\widehat{\g_{\rm free}}(v)\sum_{i=1}^d
a_i\gamma(v_i) -
\widehat{\g_{\rm free}}(w)\sum_{i=1}^d a_i\gamma(w_i)\bigg|^2 \nonumber\\
&&+ \frac{3}{2}\beta_\rho^2\int_{([0,1]^d)^2}
\frac{|\alpha_d(v,w)\widehat{\g_{\rm free}}(v,w)
+\alpha_d(v,w)\widehat{\g_{\rm free}}(v) +
\alpha_d(w,v)\widehat{\g_{\rm free}}(w)|^2}{\lambda +
\theta_d(v)+\theta_d(w)}dvdw\nonumber\\
&&+ \frac{3}{2}\int_{([0,1]^d)^2}\frac{|\delta_1(v,w)|^2}{\lambda +
\theta_d(v)+\theta_d(w)} dvdw.\nonumber
\end{eqnarray}
Note that the infimum on the six lines of (\ref{lastterm}) over
local $\g\in \bar{\C}_{d,1}$ is the same as if over $L^2$ functions
in $\bar{\C}_{d,1}$.

The strategy now follows three steps. In Step 1, we bound uniformly
in $\lambda>0$, \be \label{infoverbarC}\inf_{\g } \int_{[0,1]^d}
\frac{|-\beta_\rho-\rho \big[\sum_{i=1}^d
a_i\gamma(v_i)\big]\widehat{\g_{\rm free}}(v)|^2}{\lambda +
\theta_d(v)} + (\lambda +\theta_d(v))|\widehat{\g_{\rm free}}(v)|^2\
dv, \ee and find the $L^2$ minimizer function $\g_\lambda$.

In Step 2 we show $\g_\lambda$ is a real function and
$(\g_\lambda)_{\rm free}(0) = \sum_{|z|=1}(\g_\lambda)_{\rm free}(z)
= 0$.  Also,
we show for $x\in \Z^d$ that $\sup_{\lambda>0} |(\g_\lambda)_{\rm
free}(x)| <\infty$. Then, as $$\sup_{\lambda>0}\sup_{v\in
[0,1]^d}\frac{|\delta_0(v) +\delta_2(v)|^2}{\lambda +\theta_d(v)} \
<\ \infty \ \ {\rm and \ \ } \sup_{\lambda>0} \sup_{v,w\in [0,1]^d}
\frac{|\delta_1(v,w)|^2}{\lambda + \theta_d(v) + \theta_d(w)} \ <\
\infty,$$
 the integrals in the second and sixth lines of (\ref{lastterm}) are uniformly bounded.  Also,
 the other term in absolute value in the second line of (\ref{lastterm}), with $\g = \g_\lambda$, vanishes.

Finally, in Step 3 we show that the two integrals, with $\g =
\g_\lambda$, in the third through fifth lines of (\ref{lastterm})
are uniformly bounded in $\lambda>0$.
Hence, $\sigma^2_t(f)/t$ is uniformly bounded
over $t\geq 1$, completing the proof of Theorem \ref{maintheorem1}.
\qed

We now argue these steps. \vskip .1cm

 {\it Step 1.} By
straightforward optimizations on the quadratic expression in the
integrand, observe infimum (\ref{infoverbarC}) evaluates to
\begin{eqnarray}
\label{preliminfimum} &&{\beta_\rho^2}\int_{[0,1]^d}\frac{\lambda
+\theta_d(v)}{\rho^2\big|\sum_{i=1}^d a_i\gamma(v_i)\big |^2 +
(\lambda +\theta_d(v))^2}dv
\end{eqnarray}
with minimizer
$$\widehat{(\g_\lambda)_{\rm free}}(v) \ = \
\frac{\beta_\rho\rho \sum_{i=1}^d a_i\gamma(v_i)}
{\rho^2\big|\sum_{i=1}^d a_i\gamma(v_i)\big|^2 +(\lambda
+\theta_d(v))^2}\ .$$

We now check (\ref{preliminfimum}) is uniformly finite in
$\lambda>0$: As noted near equation (5.6) \cite{B}, which considers
almost the same integral, problems arise when $v = (0,0)$, $(0,1)$,
$(1,0)$, and $(1,1)$ in $d=2$; and in $d=1$, when $v = 0$ and $1$.

In $d=2$, by using a possible sign change, the uniform bound of
(\ref{preliminfimum}) is equivalent to bounding
$$\int_V \frac{\lambda + v_1^2 +v_2^2}{(c_1 v_1 +c_2 v_2)^2
+(\lambda +v_1^2 +v_2^2)^2}dv_1dv_2$$ where $V\in \R_+\times\R_+$ is
a neighborhood of the origin and $c_1,c_2$ are arbitrary constants
with $c_1^2 +c_2^2>0$.  As the difficulty is when $c_1v_1 +c_2v_2 =
0$, bounding the above integral is the same as bounding, with
$c_1/\sqrt{c_1^2+c_2^2} = \sin(\phi_0)$ and $c_2/\sqrt{c_1^2 +c_2^2}
= \cos(\phi_0)$,
$$\int_0^1\int_0^{\pi/2} \frac{(\lambda +r^2)r\ }{(c_1^2
+c_2^2)r^2\sin^2(\phi+\phi_0) +(\lambda +r^2)^2}d\phi dr$$ or more
simply on order
$$\int_0^1\int_0^{\pi/2} \frac{(\lambda +r^2)r\ }{(c_1^2
+c_2^2)r^2\sin^2(\phi) +(\lambda +r^2)^2}d\phi dr$$ which is finite
uniformly in $\lambda>0$ (cf. Lemma 5.2 \cite{B} for similar
calculations).

In $d=1$, (\ref{preliminfimum}) is on order $$ \int_0^1
\frac{\lambda +v^2}{v^2 + (\lambda + v^2)^2}dv$$ which also, by
straightforward computation, is finite uniformly in $\lambda>0$.

 \vskip .2cm {\it Step
2.} Noting $\overline{\gamma(r)} = -\gamma(r)$, we now show
$\widehat{(\g_\lambda)_{\rm free}}$ is the transform of a real
function:
\begin{eqnarray*}
\overline{\int_{[0,1]^d} e^{2\pi i v\cdot
x}\widehat{(\g_\lambda)_{\rm free}}(v)dv} &=&
-\int_{[0,1]^d} e^{-2\pi i v\cdot x}\widehat{(\g_\lambda)_{\rm free}}(v)dv\\
 &= & \int_{[0,1]^d}
e^{-2\pi i v\cdot x}\widehat{(\g_\lambda)_{\rm
free}}(\vec{1}-v)dv \\
& =& \int_{[0,1]^d} e^{2\pi i v\cdot x}\widehat{(\g_\lambda)_{\rm
free}}(v)dv\end{eqnarray*} where $\vec{1}$ is the vector with
components all $1$. The last sequence also shows $(\g_\lambda)_{\rm
free}$ is odd, that is $(\g_\lambda)_{\rm free}(x) =
-(\g_\lambda)_{\rm free}(-x)$ for $x\in \Z^d$. Then,
$\sum_{|z|=1}(\g_\lambda)_{\rm free}(z) = (\g_\lambda)_{\rm free}(0)
= 0$. Also, for $x\in \Z^d$, again as $(\g_\lambda)_{\rm free}$ is
odd,
\begin{eqnarray}
\label{andso} &&\sup_{\lambda>0} |(\g_\lambda)_{\rm free}(x)|\nonumber\\
&&\ \ \ \ = \ \sup_{\lambda>0} \bigg |\int_{[0,1]^d}
i\sin(2\pi v\cdot x)\widehat{(\g_\lambda)_{\rm free}}(v)dv\bigg| \\
&&\ \ \ \ \leq\  C  \sup_{\lambda>0} \int_{[0,1]^d} \frac{|\sin(2\pi
v\cdot x)||\sum_{i=1}^d a_i\gamma(v_i)|dv}{\rho^2|\sum_{i=1}^d
a_i\gamma(v_i)|^2 +(\lambda +\theta_d(v))^2}\nonumber
\end{eqnarray}
where $C=C(\rho)$. As with (\ref{preliminfimum}) above, the only
problem with the denominator in $d=2$ comes at points $v= (0,0),
(0,1),(1,0)$ and $(1,1)$, and in $d=1$ at $v=0$ and $1$.

The bound on (\ref{andso}) in $d=2$, similar to the calculation in
Step 1, is the same as bounding
$$\int_0^1\int_0^{\pi/2}
\frac{rd\phi dr}{\sin^2(\phi) +r^2}$$ which is finite. The bound on
(\ref{andso}) in $d=1$ is also finite and simpler.

\vskip .2cm {\it Step 3.} The two integrals in the third through
fifth lines of (\ref{lastterm}), after adding and subtracting
$2\beta_\rho b$ with $b= -\beta_\rho/\rho$, are bounded up to a
constant $C=C(p,\rho)$ by
\begin{eqnarray}
\label{firstintegral}
&& C\int_{([0,1]^d)^2} \frac{\big|b
-\big[\sum_{i=1}^d a_i\gamma(v_i +w_i)\big]
\widehat{(g_\lambda)_{\rm free}}(v+w)\big|^2}{\lambda +\theta_d(v) +\theta_d(w)}dvdw\\
\label{secondintegral} && +C\int_{([0,1]^d)^2}\frac{dvdw}{\lambda +
\theta_d(v)+\theta_d(w)}\\
&&\ \ \ \ \ \ \ \ \ \ \times
\bigg\{\big|b-\big[\sum_{i=1}^da_i\gamma(v_i)\big]\widehat{(g_\lambda)_{\rm
free}}(v)\big|^2 +
\big|b-\big[\sum_{i=1}^da_i\gamma(w_i)\big]\widehat{(g_\lambda)_{\rm
free}}(w)\big|^2 \bigg\} \nonumber\\
\label{thirdintegral}&&
+C\int_{([0,1]^d)^2}\frac{|\alpha_d(v,w)\widehat{(g_\lambda)_{\rm
free}}(v+w)|^2}{\lambda + \theta_d(v) +\theta_d(w)}dvdw\\
\label{fourthintegral} && +C\int_{([0,1]^d)^2}
\frac{|\alpha_d(v,w)\widehat{(g_\lambda)_{\rm free}}(v)|^2 +
|\alpha_d(w,v)\widehat{(g_\lambda)_{\rm free}}(w)|^2}{\lambda
+\theta_d(v) +\theta_d(w)}dvdw
\end{eqnarray}

The first integral (\ref{firstintegral}), noting $\big[\sum_{i=1}^d
a_i\gamma(r_i)\big]^2 = -\big|\sum_{i=1}^d a_i\gamma(r_i)\big|^2$,
is on order
$$\int_{([0,1]^d)^2} \frac{(\lambda +
\theta_d(v+w))^2}{\big|\sum_{i=1}^d a_i\gamma(v_i+w_i) \big|^2
+(\lambda + \theta_d(v+w))^2}\frac{dvdw}{\lambda
+\theta_d(v)+\theta_d(w)}$$ which in $d=2$ is bounded simply and
uniformly in $\lambda>0$ by
$$ \int_{([0,1]^2)^2} \frac{dvdw}{\theta_2(v)+\theta_2(w)} \ <
\infty .$$ In $d=1$, as $\sup_{\lambda >0}\sup_{v,w\in [0,1]}
(\lambda +\theta_1(v+w))/(\lambda +\theta_1(v)+\theta_1(w))
<\infty$, we bound on order by
$$\int_{[0,1]^2}\frac{\lambda + \theta_1(v+w)}{|\gamma(v+w)|^2
+ (\lambda + \theta_1(v+w))^2}dvdw.$$ Then, as
$$\sup_{\lambda >0}\sup_{v,w\in
[0,1]}\frac{\theta_1(v+w)}{|\gamma(v+w)|^2 +(\lambda +
\theta_1(v+w))^2} \  <\  \infty,$$ we need only bound
$$\int_{[0,1]^2} \frac{\lambda dvdw}{|\gamma(v+w)|^2 + (\lambda + \theta_1(v+w))^2} \ \leq \
\int_0^2 \frac{\lambda ds}{\sin^2(2\pi s) +(\lambda +\sin^2(\pi
s))^2}$$ which is uniformly finite in $\lambda>0$.

The second integral (\ref{secondintegral}) is analogously, and more
simply, bounded in $d=1,2$.

For the third integral (\ref{thirdintegral}), on order we need to
bound
$$\int_{([0,1]^d)^2}\frac{|\alpha_d(v,w)|^2|\sum_{i=1}^d
a_i\gamma(v_i+w_i)|^2}{[|\sum_{i=1}^d a_i\gamma(v_i+w_i)|^2 +
(\lambda + \theta_d(v+w))^2]^2}\frac{dvdw}{\lambda +
\theta_d(v)+\theta_d(w)}.$$ In $d=2$, noting the form of
$\alpha_d(v,w)$, the integral is bounded on order by
\begin{eqnarray*}
&&\int_{([0,1]^2)^2}\frac{|\sum_{i=1}^2 a_i
(\gamma(v_i)+\gamma(w_i))|^2|\sum_{i=1}^2
a_i\gamma(v_i+w_i)|^2dvdw}{[|\sum_{i=1}^2 a_i\gamma(v_i+w_i)|^2 +
\theta^2_2(v+w)]^2[\theta_2(v)+\theta_2(w)]}\\
&&\ \ \ \ \ \ \ \ \ \ \ \ \ \ \ +\int_{([0,1]^2)^2}\frac{dvdw}{\theta_2(v)+\theta_2(w)}.
\end{eqnarray*} The first
term is considered and bounded, modulo constants, in Lemma 5.3
\cite{B} through an analysis of singularities of the denominator.
The second term is clearly bounded. In $d=1$, write, for $v,w\in
[0,1]$,
\begin{eqnarray*}
\alpha_1(v,w)&=& 2 i a_1\bigg[\sin(2\pi v)[1-\cos(2\pi
 w)] +\sin(2\pi w)[1-\cos(2\pi  v)]\bigg]\\
&=& 8ia_1\sin(\pi v)\sin(\pi w)\sin(\pi(v+w)).\end{eqnarray*} Then,
the uniform bound on (\ref{thirdintegral}) follows from the bound on
the integrand
$$\sup_{\lambda>0}\sup_{v,w\in [0,1]}\frac{|\alpha_1(v,w)|^2|a_1\gamma(v+w)|^2}{[|a_1\gamma(v+w)|^2 + (\lambda +
\theta_1(v+w))^2]^2[\lambda + \theta_1(v)+\theta_1(w)]} \ < \ \infty
.$$

The last integral (\ref{fourthintegral}) is bounded on order by
$$ \int_{([0,1]^d)^2} \frac{|\alpha_d(v,w)|^2
}{|\sum_{i=1}^da_i\gamma(v_i)|^2 +
\theta_d(v)^2}\frac{dvdw}{\theta_d(v) +\theta_d(w)}. $$ As
$$\alpha_d(v,w) \ =\ \sum_{j=1}^d a_j\bigg[\gamma(v_j)(1-e^{2\pi i w_j})
+ \gamma(w_j)(1-e^{-2\pi i v_j})\bigg]$$ and $\sup_{r\in (0,1)^d}
(\sum_{j=1}^d|1-e^{\pm 2\pi i r_j}|^2)/\theta_d(r)<\infty$, the last
integral is on order \begin{equation}\label{lastintegral}
\int_{[0,1]^d}\frac{\sum_{j=1}^d|1-e^{\pm 2\pi iv_j}|^2
}{|\sum_{i=1}^da_i\gamma(v_i)|^2 + \theta_d(v)^2}dv.\end{equation}
In $d=2$, the singularities are at $v = (0,0)$, $(0,1)$, $(1,0)$ and
$(1,1)$, and as in Steps $1,2$ the bound on (\ref{lastintegral}) is
the same as
$$\int_0^1\int_0^{\pi/2}\frac{r}{\sin^2(\phi) +r^2}d\phi dr$$
which is finite. In $d=1$, as $\sup_{r\in (0,1)} |1-e^{\pm 2\pi i
r}|^2/|\gamma(r)|^2 <\infty$, the integrand in (\ref{lastintegral})
is itself finite.

\sect{Proof of Proposition \ref{FT}}
We prove the proposition in $d=2$.  The argument in $d=1$ is
analogous, and follows in particular by choosing $a_2=0$.

To make notation simple, in the following, we will omit the brackets
for singletons $\{x\}$ and two-tuple sets $\{x,y\}$ and denote them
as $x$ and $x,y$.  Also, we will drop the suffix ``${\rm nn}$'' with
respect to operators ${\a}_{n,m} = {\a}_{{\rm nn};n,m}$.  Recall
$\{e_1,e_2\}$ denotes the standard basis in $\Z^2$.

 First, from the formulas in subsection 3.1, we compute the actions
of $\a_{1,1}$ and $\a_{1,2}$ on local one-degree functions, $\g\in
\C_{2,1}$. For $x\in \Z^2\setminus\{0\}$,
\begin{eqnarray*}
(\a^e_{1,1}\g)(x) & = & (1-2\rho)\sum_{y\neq x,0}
[\g(y)-\g(x)] a_{\rm nn}(y-x) \ \ \ \ \ \ {\rm and}\\
(\a^t_{1,1}\g)(x) &= & -(1-\rho)\sum_{y\neq x,0}[\g(y)-\g(x)]a_{\rm
nn}(y-x) - \rho[\g(x)-\g(-x)]a(x)
\end{eqnarray*}
which together give
$$(\a_{1,1}\g)(x) \ =\ -\rho\sum_{y\neq x,0}
[\g(y)-\g(x)]a_{\rm nn}(y-x) -\rho[\g(x)-\g(-x)]a_{\rm nn}(x).$$
Also, for distinct $x,y\in \Z^d\setminus\{0\}$,
\begin{eqnarray*}
(\a^e_{1,2}\g)(x,y) & = & 2\beta_\rho [\g(x)-\g(y)]a_{\rm nn}(y-x)\\
(\a^t_{1,2}\g)(x,y)&= & \beta_\rho[\g(x)-\g(x-y)]a_{\rm nn}(y)+
\beta_\rho[\g(y)-\g(y-x)]a_{\rm nn}(x).\end{eqnarray*}

Then, we may write for $x\in \Z^2$ and local $\g\in \bar{\C}_{2,1}$
that $(\bar{\a}_{1,1}\g)(x)$ (cf. subsection 3.4) equals
\begin{eqnarray*}
\left\{\begin{array}{ll}
-\rho[\g(x+e_1)-\g(x-e_1)]a_1 -\rho[\g(x+e_2)-\g(x-e_2)]a_2& {\rm for \ }x\neq \pm e_1,\pm e_2,0\\
\mp\rho[\g(\pm 2e_1)-\g(\mp e_1)]a_1 -\rho[\g(e_2\pm e_1)-\g(-e_2\pm
e_1)]a_2& {\rm for \ }x= \pm e_1\\
\mp \rho[\g(\pm 2e_2)-\g(\mp e_2)]a_2  -\rho[\g(e_1\pm
e_2)-\g(-e_1\pm e_2)]a_1& {\rm for \ }x = \pm
e_2\\
0& {\rm otherwise. }\end{array}\right.
\end{eqnarray*}
 Also, for $x,y\in \Z^2$, we write (noting
 $(\bar{\a}^\cdot_{1,2}\g)(x,y)
=(\bar{\a}^\cdot_{1,2}\g)(\{x,y\})=(\bar{\a}^\cdot_{1,2}\g)(y,x)$),
\begin{eqnarray*}
(\bar{\a}^e_{1,2}\g)(x,y) &=& \left\{\begin{array}{ll}
2\beta_\rho[\g(x)-\g(x+e_1)]a_1& {\rm for \ }y=x+e_1, \ x\neq 0,-e_1\\
2\beta_\rho[\g(x)-\g(x+e_2)]a_2& {\rm for \ }y=x+e_2, \ x \neq 0,-e_2\\
0& {\rm otherwise.}
\end{array}
\right .
\end{eqnarray*}
and $(\bar{\a}^t_{1,2}\g)(x,y)$ equals\begin{eqnarray*}\left
\{\begin{array}{ll} \pm\beta_\rho[\g(x)-\g(x\mp e_1)]a_1 & {\rm for
\ }x\neq
\pm e_1, \pm e_2,0,\ y=\pm e_1\\
\pm\beta_\rho[\g(x)-\g(x\mp e_2)]a_2& {\rm for \ }x\neq
\pm e_1,\pm e_2,0, \ y = \pm e_2\\
\beta_\rho[\pm(\g(e_1) -\g(e_1\mp e_2))a_2 &\\
\ \ \ \ \ \ \ \ \ \ \ + (\g(\pm e_2) -
\g(\pm e_2 - e_1))a_1]& {\rm for \ }x= e_1, y=\pm e_2\\
\beta_\rho[\pm(\g(-e_1) -\g(-e_1\mp e_2))a_2 &\\
\ \ \ \ \ \ \ \ \ \ \ - (\g(\pm e_2) - \g(\pm e_2 +
e_1))a_1]& {\rm for \ }x= -e_1, y=\pm e_2\\
\beta_\rho[-(\g(e_1) - \g(2e_1))a_1 +(\g(-e_1) -
\g(-2e_1))a_1]& {\rm for \ }x=e_1, y=-e_1\\
\beta_\rho[-(\g(e_2) - \g(2e_2))a_2 + (\g(-e_2) -
\g(-2e_2))a_2]&{\rm for \ } x=e_2, y=-e_2\\
0&{\rm otherwise}.
\end{array}\right.
\end{eqnarray*}

We now compute corresponding Fourier transforms. To simplify
notation, we drop the subscript ``free'' and call $\g_{\rm free} =
\g$. First, we have $\widehat{(\bar{\a}_{1,1}\g)_{\rm free}}(v)$
(cf. subsection 3.5) equals
\begin{eqnarray*}
&&\sum_{x\in \Z^2} e^{2\pi i x\cdot v}
(\bar{\a}_{1,1}\g)_{\rm free}(x)\\
&& = \sum_{x\neq \pm e_1, \pm e_2, 0}-\rho e^{2\pi i x\cdot v}
[(\g(x+e_1) -
\g(x-e_1))a_1 + (\g(x+e_2) - \g(x-e_2))a_2]\\
&&\ \ \ \ \ \   -\rho e^{2\pi i v_1}[(\g(2e_1)-\g(-e_1))a_1 +
(\g(e_2+e_1) -
\g(-e_2+e_1))a_2]\\
&& \ \ \ \ \ \  -\rho e^{-2\pi i v_1}[-(\g(-2e_1)-\g(e_1))a_1
+(\g(e_2-e_1)
-\g(-e_2-e_1))a_2]\\
&&\ \ \ \ \ \ -\rho e^{2\pi i v_2} [(\g(2e_2) - \g(-e_2))a_2 +
(\g(e_1+e_2) -
\g(-e_1+e_2))a_1]\\
&&\ \ \ \ \ \ - \rho e^{-2\pi i v_2}[-(\g(-2e_2) - \g(e_2))a_2
+(\g(e_1-e_2) - \g(-e_1-e_2))a_1].
\end{eqnarray*}
The sum further equals
\begin{eqnarray*}
&& -\sum_{x\neq 0,2e_1,\atop \ \ \ \ \ \ e_1 \pm e_2,e_1}\rho
e^{2\pi i x \cdot v}e^{-2\pi i
  v_1}
\g(x) a_1 \ +\sum_{x\neq 0,-2e_1,\atop \ \ \ \ \ \  -e_1\pm e_2,
-e_1} \rho e^{2\pi i x\cdot v}
e^{2\pi i v_1}\g(x)a_1\\
&&\ \ \ \   - \sum_{x\neq 0, 2e_2,\atop \ \ \ \ \ \  e_2 \pm
e_1, e_2} \rho e^{2\pi i x\cdot v}e^{-2\pi i v_2} \g(x) a_2 \ +
\sum_{x\neq 0, -2e_2,\atop \ \ \ \
  \ \ -e_2 \pm
e_1, -e_2} \rho e^{2\pi i x\cdot v}e^{2\pi i v_2} \g(x)a_2.
\end{eqnarray*}

Recall now that $\gamma(r) = e^{2\pi i
r} - e^{-2\pi i r} = 2i\sin(2\pi r)$.
Combining and canceling terms
gives that
\begin{eqnarray*}
\widehat{(\bar{\a}_{1,1}\g)_{\rm free}}(v) &=& \rho[a_1\gamma(v_1) +a_2\gamma(v_2)]\widehat{\g}(v)\\
&& \ \ \ \ \ -\rho a_1(e^{-2\pi iv_1} -1) \g(e_1) +\rho a_1(e^{2\pi
i
  v_1}-1)\g(-e_1)\\
&& \ \ \ \ \ -\rho a_2 (e^{-2\pi i v_2} -1)\g(e_2) +\rho a_2
(e^{2\pi i v_2}
-1)\g(-e_2)\\
&&\ \ \ \ \ -\rho [ a_1\gamma(v_1) +a_2\gamma(v_2)]\g(0)\\
&=& \rho[a_1\gamma(v_1) +a_2\gamma(v_2)]\widehat{\g}(v) +
\delta_0(v)
\end{eqnarray*}
where $|\delta_0(v)| \leq \kappa(v)\sum _{|z|\leq 1}|\g(z)|$ and
$\kappa(v)$ is a bounded function on order $\kappa(v) = O(|v-z|)$
when $v\rightarrow z$ for $z= (0,0),(0,1),(1,0)$, and $(1,1)$.

We also compute that $\sqrt{2}\widehat{(\bar{\a}^e_{1,2}\g)_{\rm
free}}(v,w)$ equals
\begin{eqnarray*}
&&\sum_{x,y\in \Z^2} e^{2\pi i (x\cdot v
  +y\cdot w)}(\bar{\a}^e_{1,2}\g)_{\rm free}(x,y)\\
&&\ \ \  = 2\beta_\rho a_1 \sum_{z\neq 0,-e_1} e^{2\pi i z\cdot
(v+w)}(e^{2\pi
  i w_1} + e^{2\pi i v_1})[\g(z) - \g(z+e_1)]\\
&&\ \ \ \ \ \ \ \ \ \ \ +2\beta_\rho a_2 \sum_{z\neq 0, -e_2}
e^{2\pi i
  z\cdot (v+w)}(e^{2\pi i w_2} + e^{2\pi i v_2})[\g(z)-\g(z+e_2)]\\
&&\ \ \ = 2\beta_\rho a_1\sum_{z\neq 0, \pm e_1} e^{2\pi i
  z\cdot (v+w)}(\gamma(w_1) + \gamma(v_1))\g(z) \\
&&\ \ \ \ \ \ \ \ \ \ \ + 2\beta_\rho a_2\sum_{z\neq 0, \pm e_2}
e^{2\pi i
  z\cdot (v+w)}(\gamma(w_2) + \gamma(v_2))\g(z) \\
&&\ \ \ \ \ \ \ \ \ \ \ + 2\beta_\rho a_1
e^{2\pi i (v_1 +w_1)}(e^{2\pi i w_1} +e^{2\pi i v_1}) \g(e_1)\\
&&\ \ \ \ \ \ \ \ \ \ \ - 2\beta_\rho a_1
e^{-2\pi i (v_1 +w_1)}(e^{-2\pi i w_1} +e^{-2\pi i v_1}) \g(-e_1)\\
&&\ \ \ \ \ \ \ \ \ \ \ + 2\beta_\rho a_2 e^{2\pi i (v_2
+w_2)}(e^{2\pi i w_2}
  + e^{2\pi i v_2})\g(e_2)\\
&& \ \ \ \ \ \ \ \ \ \ \ - 2\beta_\rho a_2 e^{-2\pi i (v_2
+w_2)}(e^{-2\pi i w_2}
  +e^{-2\pi i v_2})\g(-e_2)\\
&&\ \ \ = 2\beta_\rho[a_1(\gamma(w_1) + \gamma(v_1)) +
  a_2(\gamma(w_2) + \gamma(v_2))]\widehat{\g}(v+w)\\
&&\ \ \ \ \ \ \ \ \ \ \ - 2\beta_\rho[a_1(\gamma(w_1) + \gamma(v_1))
+ a_2( \gamma(w_2) + \gamma(v_2))] \g(0)\\
&&\ \ \ \ \ \ \ \ \ \ \ + 2\beta_\rho a_1[(e^{2\pi i w_1} + e^{2\pi
i
  v_1})\g(e_1)
- (e^{-2\pi i w_1} + e^{-2\pi i v_1})\g(-e_1)]\\
&&\ \ \ \ \ \ \ \ \ \ \ +2\beta_\rho a_2[ (e^{2\pi i w_2} +e^{2\pi i
  v_2})\g(e_2) - (e^{-2\pi i w_2} + e^{-2\pi i v_2})\g(-e_2)].
\end{eqnarray*}
Also, we have $\sqrt{2}\widehat{(\bar{\a}^t_{1,2}\g)_{\rm
free}}(v,w)$ equals
\begin{eqnarray*}
&&\sum_{x,y\in \Z^2} e^{2\pi i (x\cdot v
  +y\cdot w)}(\bar{\a}^t_{1,2}\g)_{\rm free}(x,y)\\
&&\ \ \ = \beta_\rho a_1\sum_{z\neq 0,\pm e_1,\pm e_2}(e^{2\pi i
  z\cdot v}e^{2\pi i w_1} + e^{2\pi i z\cdot w}e^{2\pi i v_1})[\g(z)-\g(z-e_1)]\\
&&\ \ \ \ \ \ \ \ \ - \beta_\rho a_1\sum_{z\neq 0,\pm e_1,\pm
e_2}(e^{2\pi
  i z\cdot v}e^{-2\pi i w_1} +e^{2\pi i z\cdot w}e^{-2\pi i
  v_1})[\g(z)-\g(z+e_1)]\\
&&\ \ \ \ \ \ \ \ \ +\beta_\rho a_2\sum_{z\neq 0,\pm e_1,\pm e_2}
(e^{2\pi i
  z\cdot v}e^{2\pi i w_2} + e^{2\pi i z\cdot w}e^{2\pi i v_2})[\g(z)-\g(z-e_2)]\\
&&\ \ \ \ \ \ \ \ \ - \beta_\rho a_2\sum_{z\neq 0,\pm e_1,\pm
e_2}(e^{2\pi i
  z\cdot v}e^{-2\pi i w_2} + e^{2\pi i z\cdot w}e^{-2\pi i
  v_2})[\g(z)-\g(z+e_2)]\\
&&\ \ \ \ \ \ \ \ \ +\beta_\rho(e^{2\pi i v_1}e^{2\pi i w_2}
+e^{2\pi i
v_2}e^{2\pi i  w_1})\\
&&\ \ \ \ \ \ \ \ \ \ \ \ \ \ \ \ \ \ \ \times  [a_2(\g(e_1) -
\g(e_1-e_2)) +
a_1(\g(e_2)-\g(e_2-e_1))]\\
&& \ \ \ \ \ \ \ \ \ +\beta_\rho(e^{2\pi i v_1}e^{-2\pi i w_2} +
e^{-2\pi i
  v_2}e^{2\pi i w_1})\\
&&\ \ \ \ \ \ \ \ \ \ \ \ \ \ \ \ \ \ \ \times [-a_2(\g(e_1) -
\g(e_1+e_2)) +a_1(\g(-e_2)
  - \g(-e_2-e_1))]\\
&&\ \ \ \ \ \ \ \ \ +\beta_\rho(e^{-2\pi i v_1}e^{2\pi i w_2}+
e^{2\pi i
v_2}e^{-2\pi i  w_1})\\
&&\ \ \ \ \ \ \ \ \ \ \ \ \ \ \ \ \ \ \ \times [a_2(\g(-e_1) -
\g(-e_1-e_2))
  -a_1(\g(e_2)-\g(e_2+e_1))]\\
&&\ \ \ \ \ \ \ \ \ +\beta_\rho(e^{-2\pi i v_1}e^{-2\pi i w_2} +
e^{-2\pi i
v_2}e^{-2\pi i  w_1})\\
&&\ \ \ \ \ \ \ \ \ \ \ \ \ \ \ \ \ \ \ \times [-a_2(\g(-e_1)-\g(-e_1+e_2))-a_1(\g(-e_2) - \g(-e_2+e_1))]\\
&&\ \ \ \ \ \ \ \ \ + \beta_\rho a_1 (e^{2\pi i v_1}e^{-2\pi i w_1}
+ e^{-2\pi
  i v_1}e^{2\pi i w_1})\\
&&\ \ \ \ \ \ \ \ \ \ \ \ \ \ \ \ \ \ \ \times [-\g(e_1)+\g(2e_1) +\g(-e_1) -\g(-2e_1)]\\
&&\ \ \ \ \ \ \ \ \ + \beta_\rho a_2(e^{2\pi i v_2}e^{-2\pi i w_2} +
  e^{-2\pi i v_2}e^{2\pi i w_2})\\
&&\ \ \ \ \ \ \ \ \ \ \ \ \ \ \ \ \ \ \ \times [-\g(e_2)+\g(2e_2)
  +\g(-e_2)-\g(-2e_2)]\\
&&\ \ \ =\beta_\rho \widehat{\g}(v)[a_1 (\gamma(w_1)-\gamma(v_1+w_1))
+a_2(\gamma(w_2)-\gamma(v_2+w_2))]\\
&&\ \ \ \ \ \ \ \ \ +\beta_\rho
\widehat{\g}(w)[a_1(\gamma(v_1)-\gamma(v_1+w_1))
+a_2(\gamma(v_2)-\gamma(v_2+w_2))]\\
&&\ \ \ \ \ \ \ \ \ + \beta_\rho
\g(0)[a_1(-\gamma(v_1)-\gamma(w_1)+2\gamma(v_1+w_1)) \\
&&\ \ \ \ \ \ \ \ \ \ \ \ \ \ \ \ \ \ \ \ \ \ \ \ \ \ \ \ \ \ +
a_2(-\gamma(v_2)-\gamma(w_2) +2\gamma(v_2+w_2))]\\
&&\ \ \ \ \ \ \ \ \ - \beta_\rho a_1\g(e_1)[e^{-2\pi i w_1}
+e^{-2\pi i v_1}+ 2e^{2\pi i (v_1+w_1)}]\\
&&\ \ \ \ \ \ \ \ \ + \beta_\rho a_1\g(-e_1)[e^{2\pi i w_1} +e^{2\pi
i v_1} +2e^{-2\pi i (v_1+w_1)}]\\
&&\ \ \ \ \ \ \ \ \ - \beta_\rho a_2\g(e_2)[e^{-2\pi i w_2} +
e^{-2\pi i v_2}+2e^{2\pi i (v_2+w_2)}]\\
&&\ \ \ \ \ \ \ \ \ + \beta_\rho a_2\g(-e_2)[e^{2\pi i w_2} +
e^{2\pi i v_2}+ 2e^{-2\pi i (v_2+w_2)}]
\end{eqnarray*}

Putting terms together, we compute the Fourier transform for
$(\bar{\a}_{1,2}\g)_{\rm free}$. In particular, we note the last
four lines of the computation for
$\widehat{(\bar{\a}^t_{1,2}\g)_{\rm free}}$ and the last two lines
of the computation for $\widehat{(\bar{\a}^e_{1,2}\g)_{\rm free}}$
are $O(1)$ as $(v,w)\rightarrow (z_1,z_2)$ for $z_1,z_2 =
(0,0),(0,1),(1,0)$ and $(1,1)$.
However, they match in the sense, when they are added to each other,
the sum is small on the desired order. We have
\begin{eqnarray*}
&&\sqrt{2}\widehat{(\bar{\a}_{1,2}\g)_{\rm free}}(v,w)\\
 &&\ \ \ \ \ =
2\beta_\rho[a_1(\gamma(w_1) + \gamma(v_1)) +
  a_2(\gamma(w_2) + \gamma(v_2))]\widehat{\g}(v+w)\\
&&\ \ \ \ \ \ \ + \beta_\rho \widehat{\g}(v)[a_1 (\gamma(w_1)-\gamma(v_1+w_1))
+a_2(\gamma(w_2)-\gamma(v_2+w_2))]\\
&&\ \ \ \ \ \ \ +\beta_\rho \widehat{\g}(w)[a_1(\gamma(v_1)-\gamma(v_1+w_1))
+a_2(\gamma(v_2)-\gamma(v_2+w_2))] +\delta_1(v,w)
\end{eqnarray*}
where $|\delta_1(v,w)| \leq \kappa(v,w)\sum_{|z|\leq 1}|\g(z)|$ and
$\kappa(v,w)$ is a bounded function on order $\kappa^2(v,w) =
O(|v-z_1|^2 + |w-z_2|^2)$ when $(v,w)\rightarrow (z_1,z_2)$ for $z_1,z_2 = (0,0),(0,1),(1,0)$ and $(1,1)$.
\vskip .2cm

{\bf Acknowledgement.}  I thank M. Loulakis for discussions on
Corollary 1 in \cite{Loulakis} during the stimulating conference in
Budapest, August 22-26, 2005, ``Large scale behavior of interacting
particle
               systems: Fluctuations and hydrodynamics,'' organized
               by C. Landim, J. Fritz, B. Toth and B. Valko.

\bibliographystyle{plain}

\end{document}